

\input epsf.tex

\def\2{{1\over 2}}

\def\d{\delta}
\def\a{\alpha}
\def\b{\beta}
\def\g{\gamma}

\def\s{\sigma}
\def\e{\epsilon}
\def\l{\lambda}
\def\o{\omega}

\def\fun#1#2#3{#1\colon #2\rightarrow #3}
\def\norm#1{\Vert #1 \Vert}

\def\frac#1#2{{{#1} \over {#2}}}

\def\st{\;\colon\;}
\def\tends{\rightarrow}

\def\dx{\hbox{{\rm d}$x$}}

\def\dr{ {\rm d} }
\def\dt{\hbox{{\rm d}$t$}}

\def\R{{\bf R}}
\def\N{{\bf N}}
\def\Z{{\bf Z}}

\def\T{{\bf T}}
\def\Space{{\bf S}}

\def\H{{\cal H}}
\def\M{{\cal M}}
\def\mod{{\rm mod}}
\def\thm#1{\vskip 1 pc\noindent{\bf Theorem #1.\quad}\sl}
\def\lem#1{\vskip 1 pc\noindent{\bf Lemma #1.\quad}\sl}
\def\prop#1{\vskip 1 pc\noindent{\bf Proposition #1.\quad}\sl}
\def\cor#1{\vskip 1 pc\noindent{\bf Corollary #1.\quad}\sl}

\def\proof{\rm\vskip 1 pc\noindent{\bf Proof.\quad}}
\def\fin{\par\hfill $\backslash\backslash\backslash$\vskip 1 pc}
\def\txt#1{\quad\hbox{#1}\quad}
\def\m{\mu}
\def\L{{\cal L}}
\def\H{{\cal H}}

\def\s{\sigma}

\def\o{\omega}
\def\r{\rho}

\def\2{\frac{1}{2}}
\def\inn#1#2{{\langle #1 ,#2\rangle}}

\def\W{{  {\cal W}  }}
\def\V{{  {\cal V}  }}

\def\Group{{  {\rm Group}  }}

\def\part{{\partial_{x}}}

\def\pprime{{{}^\prime{}^\prime}}

\def\I{{\cal I}_{(U^-,U^+)}}



\baselineskip= 17.2pt plus 0.6pt
\font\titlefont=cmr17
\centerline{\titlefont The Aubry set}
\vskip 1 pc
\centerline{\titlefont for a version of the Vlasov equation}
\vskip 4pc
\font\titlefont=cmr12
\centerline{         \titlefont {Ugo Bessi}\footnote*{{\rm 
Dipartimento di Matematica, Universit\`a\ Roma Tre, Largo S. 
Leonardo Murialdo, 00146 Roma, Italy.}}   }{}\footnote{}{
{{\tt email:} {\tt bessi@matrm3.mat.uniroma3.it} Work partially supported by the PRIN2009 grant "Critical Point Theory and Perturbative Methods for Nonlinear Differential Equations.}} 
\vskip 0.5 pc
 
\par
\vskip 2pc
\centerline{\bf Abstract}
We check that several properties of the Aubry set, first proven for  finite-dimensional Lagrangians by Mather and Fathi, continue to hold in the case of the infinitely many interacting particles of the Vlasov equation on the circle.

\vskip 2 pc
\centerline{\bf  Introduction}
\vskip 1 pc
The Vlasov equation on the circle governs the motion of a group of particles on $S^1\colon=\frac{\R}{\Z}$ under the action of an external potential $V(t,x)$ and a mutual interaction $W$. More precisely, we let $I=[0,1)$, we lift our particles to $\R$, and we  parametrize them at time $t$ by a function $\s_t\in L^2(I,\R)$; we require that 
$\s_t$ satisfies the differential equation in $L^2(I,\R)$
$$\ddot\s_tz= -V^\prime(t,\s_t z)-
\int_I W^\prime(\s_tz-\s_t\bar z)\dr\bar z  .  
\eqno (ODE)_{Lag}$$
Our standing hypotheses on the potentials $V$ and $W$ are

\vskip 1pc

\noindent $\bullet$ $V\in C^2(S^1\times S^1)$, $W\in C^2(S^1)$; moreover $W$, seen as a function on $\R$, is even; up to adding a constant, we can suppose that $W(0)=0$.

\vskip 1pc

There is an element of arbitrariness in choosing the lift of the particles to $\R$ and in parametrizing them; that's why we are less interested in the evolution of the labelling $\s_t$ than in the evolution of the measure it induces. In other words, we want to study the measures on $S^1\times\R$ given by 
$\mu_t\colon=(\pi\circ\s_t,\dot\s_t)_\sharp\nu_0$, where $\nu_0$ denotes the Lebesgue measure on $I$, $\fun{\pi}{\R}{S^1}$ is the natural projection and $(\cdot)_\sharp$ denotes the push-forward. A standard calculation shows that, if $\s_t$ satisfies 
$(ODE)_{Lag}$, then $\mu_t$ satisfies, in the weak sense,
$$\partial_t\mu_t+v\partial_x\mu_t=\partial_v(\mu_t\partial_x P_t)
\eqno (ODE)_{meas}$$
where 
$$P_t(x)=
V(t,x)+\int_{S^1\times\R}W(x-\bar x)\dr\mu_t(\bar x,v)=
V(t,x)+\int_I W(x-\s_t\bar z)\dr\bar z  .   $$

Problem $(ODE)_{meas}$ (see [10], [9]) is Lagrangian; actually, many results of Aubry-Mather theory can be extended to curves of measures which are "minimal" in a suitable way. Here, however, we follow the approach of [8]: we are going to work with 
$(ODE)_{Lag}$, keeping track of its symmetries. Quotienting 
$(ODE)_{Lag}$ by its symmetry group, we shall get a problem equivalent to $(ODE)_{meas}$. Though in this paper we restrict ourselves to the one-dimensional situation, we recall that L. Nurbekian (see [14]) has extended the results of [8] about minimal parametrizations to tori of any dimension.

The aim of this paper is to do a few simple checks, showing that many features of Aubry-Mather theory persist in this setting; actually, we shall check that the main theorems of [7] continue to hold. In section 1, we recall the main results of [8] on $(ODE)_{Lag}$ and its symmetries; in section 2, following [8], [10] and [7], we define the Hopf-Lax semigroup and we show that it has fixed points. We also show that the value function satisfies the Hamilton-Jacobi equation on $L^2$. In section 3, we show that 
$(ODE)_{Lag}$ admits invariant measures minimal in the sense of Mather; as a consequence, we can define Mather's conjugate actions $\a$ and $\b$. In section 4, we recall two different definitions of the Aubry set, one of Mather's and the other of Fathi's; we show that, also in this case, the two definitions coincide. In section 5, we shall show that $(ODE)_{meas}$ admits a solution $\mu_t$ which is periodic (i. e. $\mu_0=\mu_1$) and has irrational rotation number. We shall see that, as a consequence of the KAM theorem, if the rotation number $\o$ is sufficiently irrational, and $V$ and $W$ are sufficiently regular and small (depending on $\o$), then $\mu_t$ has a smooth density.

\vskip 2pc
\centerline{\bf \S 1}
\centerline{\bf Notation and preliminaries}

\vskip 1pc

Since $V$ and $W$ are periodic, we have that $(ODE)_{Lag}$ is invariant by the action of $L^2_\Z\colon=L^2(I,\Z)$; in other words, if $\s_t$ is a solution and $h\in L^2_\Z$, then $\s_t+h$ is a solution too. Moreover, $(ODE)_{Lag}$ is invariant by the group $G$ of the measure-preserving transformations of $I$ into itself; indeed, such maps do not change the value of the integral defining $P_t(x)$. An idea of [8] is to quotient $L^2(I)$ by these two groups; we recall from [8] some facts about this quotient. 

We shall denote by $\norm{\cdot}$ the norm on $L^2(I)$, 
$\inn{\cdot}{\cdot}$ the internal product. We set
$$\T\colon =\frac{L^2(I)}{L^2_\Z(I)} . $$
The space $\T$ is metric, with distance between the equivalence classes $[M]$ and $[\bar M]$ given by
$$dist_\Z([M],[\bar M])=
\inf_{Z\in L^2_\Z}
||
M-\bar M-Z
||=
||
|M-\bar M|_{S^1}
||$$
where 
$$|m|_{S^1}\colon=\min_{k\in\Z}|m+k| . $$ 
We note that, for each $x\in I$, we can measurably choose 
$Zx\in\Z$ such that $|Mx-\bar Mx-Zx|=|Mx-\bar Mx|_{S^1}$; as a consequence, we get that the $\inf$ in the definition of $dist_\Z$ is a minimum; we also get the second equality above.

Let $\Group$ denote the group of the measure-preserving transformations of $I$ with measurable inverse; for 
$M,\bar M\in L^2(I)$ we set
$$dist_{weak}(M,\bar M)=\inf_{G\in\Group}
dist_\Z(M\circ G,\bar M) . $$
This yields that $M$ and $M\circ G$, which we would like to consider equivalent, have zero distance; however, if we say that 
$M\simeq\bar M$ when $\bar M=M\circ G$ for some 
$G\in\Group$, then the equivalence classes are not closed in 
$\T$, essentially because the $\inf$ in the definition of 
$dist_{weak}$ is not a minimum: it is possible (see [8]) that 
$dist(M,\bar M)=0$ even if $M$ and $\bar M$ are not equivalent.  But we can consider their closure if we look at the equivalence relation from the right point of view, i. e. that of the measure induced by $M$. 

We denote by $\rm{Meas}$ the space of Borel probability measures on 
$S^1$, and we let $\fun{\pi}{\R}{S^1}$ be the natural projection. We introduce the map
$$\fun{\Phi}{L^2(I)}{{\rm Meas}}, \qquad
\fun{\Phi}{M}{(\pi\circ M)_\sharp\nu_0}$$
where $(\cdot)_\sharp$ denotes push-forward and $\nu_0$ is the Lebesgue measure on $I$. We note that $\Phi$ is invariant under the action of $L^2_\Z$ and $\Group$; in other words, if 
$Z\in L^2_\Z$ and $G\in\Group$, then 
$\Phi(u)=\Phi((u+Z)\circ G)$. We say that 
$M\simeq\bar M$ if $\Phi(M)=\Phi(\bar M)$. We set 
$\Space\colon=\frac{L^2(I)}{\simeq}$; on this space, we consider the metric
$$dist_\Space([M],[\bar M])=\inf\{
||M^\ast-\bar M^\ast||   
\st
M^\ast\in[M],\quad \bar M^\ast\in[\bar M]  
\}    .   $$
The infimum above is a minimum: one can always find a minimal couple $(M^\ast,\bar M^\ast)$ with $M$ monotone and taking values in $[0,1]$, and $\bar M$ monotone and taking values in 
$[-\frac{3}{2},\frac{3}{2}]$. Now $\Space$ is isometric to the space of Borel probability measures on $S^1$ with the 
2-Wasserstein distance; in particular, it is a compact space.

Another fact proven in [8] is that 
$dist_{weak}(M,\bar M)=dist_\Space([M],[\bar M])$.

By proposition 2.9 of [8], which we copy below, the 
$L^2_\Z$-equivariant (or $L^2_\Z$-equivariant and 
$\Group$-equivariant) closed forms on $L^2(I)$ have a particularly simple structure: the first equivariant cohomology group of 
$L^2(I)$ is $\R$.

\prop{1.1} Let $\fun{S}{L^2(I)}{\R}$ be $C^1$.

\noindent 1) If $\dr S$ is $L^2_\Z$-periodic in the sense that 
$\dr_{M+Z}S=\dr_M S$ for all $Z\in L^2_\Z(I)$, then there is a unique $C\in L^2(I)$ and a function $\fun{s}{L^2(I)}{\R}$, of class $C^1$ and $L^2_\Z$-periodic, such that
$$S(M)=s(M)+\inn{C}{M}   .  $$

\noindent 2) If, in addition, $\fun{}{M}{\dr_M S}$ is 
rearrangement-invariant (i. e. if $\dr_MS=\dr_{M\circ G}S$ for all 
$G\in\Group$), then $C$ is constant and $s$ is 
rearrangement-invariant.

\rm

\vskip 1pc

In view of the lemma above, for $c\in\R$ we define the Lagrangian 
$\L_c$ as
$$\fun{\L_c}{S^1\times L^2(I)\times L^2(I)}{\R},\qquad
\L_c(t,M,N)=\2||N||^2-\inn{c}{N}-\V(t,M)-\W(M)$$
where
$$\V(t,M)=\int_I V(t,Mx)\dx ,
\txt{and}
\W(M)=\2\int_{I\times I}W(Mx-Mx^\prime)\dx\dx^\prime . $$
In order to define the $c$-minimal orbits of $\L$, we let 
$K\subset\R$ be an interval; following [1], we say that 
$u\in L^1(K,L^2(I))$ is absolutely continuous if there is 
$\dot u\in L^1(K,L^2(I))$ such that, for any 
$\phi\in C^1_0(K,\R)$, we have that
$$\int_K u_t(x)\dot\phi(t)\dt=-\int_K\dot u_t(x)\phi(t)\dt . 
\eqno (1.1)$$
The equality above is an equality in $L^2(I)$, i. e. it holds for a. e. 
$x\in I$; however, it is easy to see that the exceptional set does not depend on $\phi$, and thus that, for a. e. $x$, the map 
$\fun{}{t}{u_t(x)}$ is A. C. with derivative $\dot u_t(x)$. We shall denote by $AC(K,L^2(I))$ the class of A. C. functions from $K$ to 
$L^2(I)$.

Let $c\in\R$; we say that $\s\in AC(K,L^2(I))$ is $c$-minimal for 
$\L$ if, for any interval $[t_0,t_1]\subset K$ and any
$\tilde\s\in AC((t_0,t_1),L^2(I))$ satisfying 
$$\tilde\s_{t_1}-\s_{t_1}\in L^2_\Z(I) \txt{and} 
\tilde\s_{t_2}-\s_{t_2}\in L^2_\Z(I),$$ 
we have that
$$\int_{t_0}^{t_1}\L_c(t,\s_t,\dot\s_t)\dt\le
\int_{t_0}^{t_1}\L_c(t,\tilde\s_t,\dot{\tilde\s}_t)\dt  .  $$
We forego the standard proof that $c$-minimal orbits solve 
$(ODE)_{Lag}$.

Let now $n\in\N$, and let ${\cal A}_n$ be the $\s$-algebra on $I$  generated by the intervals $[\frac{i}{n},\frac{i+1}{n})$ with
$i\in(0,\dots,n-1)$; we call ${\cal C}_n$ the closed subspace of the 
${\cal A}_n$-measurable functions of $L^2(I)$, and we denote by 
$\fun{P_n}{L^2(I)}{{\cal C}_n}$ the orthogonal projection.  We have a bijection
$$\fun{D_n}{\R^n}{{\cal C}_n},\qquad
\fun{D_n}{(q_1,\dots,q_n)}{
\sum_{i=0}^{n-1}q_i1_{[\frac{i}{n},\frac{i+1}{n})}(x)
}    .   $$

We also note that the space 
$S^1\times{\cal C}_n\times{\cal C}_n$ is invariant by the Euler-Lagrange flow of $(ODE)_{Lag}$.

\vskip 2pc

\centerline{\bf \S 2}
\centerline{\bf The Hopf-Lax semigroup}

\vskip 1pc

\noindent{\bf Definitions.} Let us denote by $C_{\Group}(\T)$ the set of functions 
$U\in C(L^2(I),\R)$ which are $L^2_\Z$ and $\Group$ equivariant. It is standard (proposition 2.8 of [8]) each 
$U\in C_{\Group}(\T)$ quotients to a continuous function on the compact space ${\bf S}$; in particular, it is bounded. 

Given $M\in L^2(I)$, $U\in C_{\Group}(\T)$ and $t>0$, we define
$$(A^t_c U)(M)=
\inf\{
\int_{0}^t\L_c(s,\s_s,\dot\s_s)\dr s+U(\s_{0})\st 
\s\in AC([0,t],L^2(I)),\quad\s_t=M
\}   .   \eqno (2.1)  $$

We shall denote by $Mon$ the space of the maps 
$\fun{\s}{I}{\R}$ which are monotone increasing and satisfy 
$\s(1-)\le\s(0)+1$. We endow $Mon$ with the topology it inherits from $L^2(I)$, which turns it into a locally compact space.

We group together the statements of a few lemmas of [8] and [10]; for a slightly different proof, point 1) is lemma 2.1 of [5], point 2) is lemma 2.8, point 4) proposition 2.2.

\prop{2.1} Let $U\in C_{\Group}(\T)$, let $t>0$ and let 
$\fun{A^t_cU}{L^2(I)}{\R}$ be defined by (2.1); then, the following statements hold.

\noindent 1) $A^t_c U$ is $L^2_\Z$ and 
$\Group$-equivariant.

\noindent 2) $A^t_c U$ is $L(t)$-Lipschitz for $dist_{weak}$ (or for $dist_{\bf S}$, since we have seen that the two distances coincide). The constant $L(t)$ does not depend on $U$. Moreover, $L(t)\le L$ for $t\ge 1$.

\noindent 3) As a consequence of 1) and 2), 
$A^t_c U\in C_{\Group}(\T)$.

\noindent 4) Let $M\in Mon$; then, the $\inf$ in (2.1) is a minimum; more precisely, there is $\s\in AC([0,t],L^2(I))$ with $\s_t=M$, 
$\s_s\in Mon$ for $s\in[0,t]$ and such that
$$(A^t_c U)(M)=\int_{0}^t\L_c(s,\s_s,\dot\s_s)\dr s+U(\s_{0}) . $$
The function $\s$ is $c$-minimal on $({0},t)$ and solves 
$(ODE)_{Lag}$.

\noindent 5) Since $\L_c$ is one-periodic in time, $A^t_c$ has the semigroup property on the integers: in other words, if $t>0$ and 
$s\in\N$, then
$$A^{t+s}_cU=A^t_c(A^s_c U)  .  $$

\rm
\vskip 1pc

Let $\l\in\R$; by point 3) of the last lemma, we can define a map
$$\fun{\Lambda_{c,\l}}{C_{\Group}(\T)}{C_{\Group}(\T)}$$
$$\fun{\Lambda_{c,\l}}{U}{(A^1_c U)(\cdot)+\l}  .  $$
It follows immediately from the definition of $A^1_c U$ that

\vskip 1pc

\noindent $\bullet$ $\Lambda_{c,\l}$ is monotone, i. e., if 
$U_1\le U_2$, then $\Lambda_{c,\l}U_1\le\Lambda_{c,\l}U_1$.

\noindent $\bullet$ If $a\in\R$, then 
$\Lambda_{c,\l}(U+a)=\Lambda_{c,\l}U +a$.

\vskip 1pc

\noindent These two facts easily imply that

\vskip 1pc

\noindent $\bullet$ $\Lambda_{c,\l}$ is continuous (actually, 
$1$-Lipschitz) from 
$C_\Group(\T)$ to itself, if we put on $C_\Group(\T)$ the $\sup$ norm. 

\vskip 1pc

Again, we refer the reader to [8], [10] (or to [7], since the finite dimensional proof is the same) for the next lemma; in [5], point 1) is proposition 2.11. Point 2) follows in a standard way by point 1) and the semigroup property.

\prop{2.2} 1) There is a unique $\l\in\R$ (which we shall call $\a(c)$) such that $\Lambda_{c,\l}$ has a fixed point in 
$C_{\Group}(\T)$. By point 2) of proposition 2.1, any fixed point is 
$L$-Lipschitz.

\noindent 2) Let $U$ be a fixed point of $\Lambda_{c,\l}$, and let 
$M\in Mon$. Then, there is $\s\in AC_{loc}((-\infty,0],L^2(I))$ such that $\s_t\in Mon$ for $t\in(-\infty,0)$, $\s_0=M$ and, for all 
$k\in\N$, 
$$U(M)=\int_{-k}^0[\L_c(t,\s_t,\dot\s_t)+\a(c)]\dt+U(\s_{-k}) . $$
The function $\s$ is $c$-minimal on $(-\infty,0)$ and solves 
$(ODE)_{Lag}$.

\rm
\vskip 1pc

Now we introduce the notation of [7] for the 
Hopf-Lax semigroups, forward ($T^-_{t}$) and backward 
($T^+_{-t}$) in time. The signs $+$ and $-$ point, apparently, in the wrong direction; a possible justification is that, when the semigroup goes forward in time, the characteristics go backward, and vice-versa.

\vskip 1pc

\noindent {\bf Definition.} Let $U\in C_{\Group}(\T)$, let 
$M\in L^2(I)$ and let $\a(c)$ be as in proposition 2.2; for $t\ge 0$, we define
$$(T^-_{t}U)(M)=\inf\{
U(\g_{0})+\int^{t}_0[
\L_c(s,\g_s,\dot\g_s)+\a(c)
]\dr s  \st \g_t=M
\}  $$
and 
$$(T^+_{-t}U)(M)=\sup\{
U(\g_{0})-\int_{-t}^0[
\L_c(s,\g_s,\dot\g_s)+\a(c)
]\dr s  \st \g_{-t}=M
\}  .   $$

\vskip 1pc

We note that, by proposition 2.1, $T^-_tU$ and $T^+_{-t}U$ belong to $C_\Group(\T)$. By proposition 2.2, 
$T^-_1=\Lambda_{c,\a(c)}$ has a fixed point; we cannot say the same for $T^+_{-1}$ because the choice $\l=\a(c)$, which yields a fixed point of $T^-_1$, may not yield a fixed point of $T^+_{-1}$; we shall have to wait until theorem 4.2 below to see that this is actually the case, and that both operators have fixed points.

By point 5) of proposition 2.1, if $U$ is a fixed point of $T_1^-$, then, for $t\ge 0$, $T^-_{t+1}U=T_tU$; in other words, the function $(T^-_{t}U)(M)$ defined on 
$[0,+\infty)\times L^2(I)$ can be extended by periodicity to 
$\R\times L^2(I)$. 
As a final remark, if $M\in Mon$, it follows by proposition 2.1 that 
$(T^-_tU)(M)$ and $(T^+_{-t}U)(M)$ are a minimum and a maximum respectively.

\vskip 1pc
\noindent{\bf Definition.} We shall say that a function 
$U\in C_{\Group}(\T)$ is $c$-dominated if, for every $m<n\in\Z$ and every $\s\in AC([m,n],L^2(I))$, we have that
$$U(\s_n)-U(\s_m)\le
\int_m^n[
\L_c(t,\s_t,\dot\s_t)+\a(c)
]\dt  .  $$
We note that there are $c$-dominated functions: for instance, the fixed points of $T^-_1$, given by proposition 2.2, are 
$c$-dominated by formula 2.1.

\vskip 1pc

\noindent{\bf Definition.} If $\s\in AC([a,b],L^2(I))$ and 
$\s_t\in Mon$ for $t\in[a,b]$, we shall say that 
$\s\in AC_{mon}([a,b])$. By point 4) of proposition 2.1, if 
$M\in Mon$ there is $\s\in AC_{mon}$ minimal (or maximal) in the definition of $T^-_{t}U(M)$ (or of $T^+_{-t}U(M)$.)

\vskip 1pc

\noindent{\bf Definition.} Let $U\in C_\Group(\T)$ be 
$c$-dominated and let $a<b\in\Z\cup\{ \pm\infty \}$; we say that 
$\g\in AC_{mon}([a,b])$ is calibrating if, for any 
$[m,n]\subset[a,b]$ with $m$ and $n$ integers, we have
$$U(\g_n)-U(\g_m)=\int_m^n[
\L_c(t,\g_t,\dot\g_t)+\a(c)
]\dt   .   $$
It follows from (2.1) that a calibrating function $\g$ is $c$-minimal on $[a,b]$, and thus it satisfies $(ODE)_{Lag}$.

\vskip 1pc

We state at once a relation between these definitions; it comes, naturally, from [7].

\lem{2.3} 1) Let $U\in C_{\Group}(\T)$; then $U$ is $c$-dominated iff $U\le T^-_{n}U$ (or iff $T^+_{-n}U\le U$) for all $n\ge 0$.

\noindent 2) Moreover, $T^-_{n}(U)=U$ (or $T^+_{-n}U=U$) for all $n\in\N$ iff 
$U$ is $c$-dominated and, for each $M\in Mon$, there is a calibrating curve $\g\in AC_{mon}((-\infty,0])$ (or 
$\g\in AC_{mon}([0,+\infty))$ with $\g_0=M$.

\proof Point 1) is a rewording of the definition of $c$-dominated. We prove point 2); if $T^-_{n}U=U$, then $U$ is $c$-dominated by point 1); the existence of a calibrating curve $\g$ follows from point 2) of proposition 2.2. To prove the converse, let $M\in Mon$ and let 
$\g$ be calibrating on $(-\infty,0]$ with $\g_0=M$; then, 
$$U(M)-U(\g_{-1})=U(\g_0)-U(\g_{-1})=\int_{-1}^0[
\L_c(t,\g_t,\dot\g_t)+\a(c)
]\dt  .  $$
By the definition of $T_{1}^-$, this means that 
$(T_{1}^-U)(M)\le U(M)$; since the opposite inequality holds by point 1), we have that $(T^-_{1}U)(M)=U(M)$ for all $M\in Mon$. Since $U$ is continuous and equivariant, and since by [8] any 
$N\in L^2(I)$ can be approximated by $M\circ G_n+Z_n$ with 
$M\in Mon$, $G_n\in\Group$ and $Z_n\in L^2_\Z$, we have that 
$(T^-_1U)(N)=U(N)$ for all $N\in L^2(I)$, and we are done.

\fin

\vskip 1pc

The Lagrangian $\L_c$ has a Legendre transform $\H_c$; an easy calculation shows that
$$\fun{\H_c}{S^1\times L^2(I)\times L^2(I)}{\R}$$
$$\H_c(t,\s,p)=\2||c+p||_{L^2(I)}^2+\V(t,\s)+\W(\s) . $$

What we really need are subsolutions of Hamilton-Jacobi; that's why we give the following definition.

\vskip 1pc

\noindent{\bf Definition.} We define $Mon_3$ as the set of monotone functions $\g$ on $I$ such that $\g(1-)\le\g(0)+3$. Let 
$\fun{U}{\R\times L^2(I)}{\R}$, and let $M\in Mon$. We say that 
$(a,\xi)\in\R\times L^2(I)$ is the "lazy differential" of $U$ at 
$(t,M)$ if there is $K>0$ such that
$$U(t+h,M+N)-U(t,M)\le ah+\inn{\xi}{N}+
K(|h|^2+||N||^2)\qquad
\forall(h,N)\in\R\times L^2(I)  \leqno i)$$
and
$$U(t+h,M+N)-U(t,M)\ge ah+\inn{\xi}{N}+o(|h|+||N||)
\leqno ii)$$
for all $N\in L^2(I)$ such that $M+N\in Mon_3$.

We set $(\partial_t U(t,M),\partial_M U(t,M))\colon=(a,\xi)$.

\lem{2.4} If $U$ is lazily differentiable at $(t,M)$, then the lazy differential $(a,\xi)$ is unique.

\proof Let $(a^\prime,\xi^\prime)$ be another lazy differential; if we set $N=0$ in $i)$, $ii)$, we get that $a=a^\prime$. 

If we set $h=0$, $N=\e\bar N$ for $\e>0$, and we subtract $ii$) 
$$U(t,M+\e\bar N)-U(t,M)\ge\e\inn{\xi^\prime}{\bar N}+ o(\e)$$
from $i$)
$$U(t,M+\e\bar N)-U(t,M)\le\e\inn{\xi}{\bar N}+K\e^2  ,  $$
we get that
$$\inn{\xi-\xi^\prime}{\bar N}\ge 0$$
for all $\bar N$ such that $M+\e\bar N\in Mon_3$ for $\e$ positive and small. Exchanging the r\^oles of $\xi$ and $\xi^\prime$, we get that
$$\inn{\xi-\xi^\prime}{\bar N}=0$$
for all $\bar N$ such that $M+\e\bar N\in Mon_3$ for $\e$ positive and small. In particular, the formula above holds for 
$\bar N=1_{[c,1]}$ and $\bar N_1=1_{[d,1]}$; subtracting, we get that
$$\int_c^d(\xi(x)-\xi^\prime(x))\dr\nu_0(x)=0\qquad
\forall 0\le c< d\le 1  .  $$
Thus, $\xi=\xi^\prime$, as we wanted.

\fin

\prop{2.5} Let $U\in C_\Group(\T)$. For $t>0$, let us set 
$\hat U(t,M)=(T^-_tU)(M)$. For $(t,M)\in(0,+\infty)\times Mon$, let us suppose that there is a unique curve $\s$ such that $\s_t=M$ and 
$$\hat U(t,\s_t)-U(\s_0)=\int_0^t[
\L_c(s,\s_s,\dot\s_s)+\a(c)
]\dr s   .  \eqno (2.2)$$
Then, $\hat U$ is lazily differentiable at $(t,M)$ and 
$$\partial_t\hat U(t,M)+\H_0(t,M,c+\partial_M\hat U(t,M))=\a(c) . 
\eqno (2.3)$$
As a partial converse, if $\hat U$ is Fr\'echet differentiable at 
$(t,M)\in(0,+\infty)\times Mon$, then there is a unique $\s$ minimal in (2.2), which satisfies (2.3) by the statement above.

\proof The proof is identical to the finite-dimensional one. We begin with the converse. 

Let $\hat U$ be Fr\'echet differentiable at 
$(t,M)\in(0,+\infty)\times Mon$; by proposition 2.1, there is a curve $\s$ such that (2.2) holds; we want to prove that it is unique. For $N\in L^2(I)$, let us set
$$\tilde\s_s=\s_s+(N-M)\frac{s}{t}  .  $$
Since $\tilde\s_t=N$, $\tilde\s_0=\s_0$ and $\s$ is minimal, the definition of $\hat U$ implies the first inequality below.
$$\hat U(t,N)-\hat U(t,M)\le
\int_0^t[\L_c(s,\tilde\s_s,\dot{\tilde\s}_s)-\L(s,\s_s,\dot\s_s)]\dr s
\le$$
$$\int_0^t[
\inn{\dot\s_s-c}{\frac{N-M}{t}}-
\inn{\V^\prime(s,\s_s)+\W^\prime(\s_s)}{\frac{(N-M)s}{t}}
]\dr s
+K||N-M||^2 = $$
$$\inn{\dot\s_t-c}{N-M}+K||N-M||^2  .  \eqno (2.4)$$
The second inequality above comes from a Taylor development of 
$\L_c$, and from the fact that the second derivatives of $V$ and $W$ are bounded; the equality comes from an integration by parts and the fact that $\s$, by point 4) of proposition 2.1, solves 
$(ODE)_{Lag}$.

If $\hat U$ is Fr\'echet differentiable at $(t,M)$, the last formula implies that
$$\partial_M\hat U(t,M)=(\dot\s_t-c)  .  \eqno (2.5)$$
Since $\s$ satisfies (2.2), it is calibrating, and thus it solves 
$(ODE)_{Lag}$; we have just seen that its final speed at $t$ satisfies the formula above; since the existence and uniqueness theorem holds for $(ODE)_{Lag}$, we get that the minimizer at 
$(t,M)$ is unique. It remains to prove that (2.3) holds; since we have just shown that the minimizer $\s$ is unique, this follows from the direct statement, which we presently prove.

Let us suppose that $(t,M)\in(0,+\infty)\times Mon$, and let the minimum in (2.2) be attained on a unique $\s$. We want to prove that $\hat U$ is lazily differentiable and satisfies (2.3) at 
$(t,M)$.
For $h\in\R$ and $N\in L^2(I)$, we set
$$\tilde\s_s=\s_s+\frac{N-M}{t+h}s
-(\s_{t+h}-\s_t)\frac{s}{t+h}$$
and we see that $\tilde\s_{t+h}=N$ while 
$\tilde\s_0=\s_0$. We get as above that
$$\hat U(t+h,N)-\hat U(t,M)\le
\int^{t+h}_0[\L_c(s,\tilde\s_s,\dot{\tilde\s}_s)+\a(c)]\dr s-
\int^{t}_0[\L_c(s,\s_s,\dot\s_s)+\a(c)]\dr s= $$
$$\int_t^{t+h}[
\L_c(s,\s_s,\dot\s_s)+\a(c)
]\dr s+
\int_0^{t+h}[
\L_c(s,\tilde\s_s,\dot{\tilde\s}_s)-\L_c(s,\s_s,\dot\s_s)
]\dr s  .  $$
We also note that, since 
$||V^\prime||_\infty+||W^\prime||_\infty\le K$, we have
$||\V^\prime||+||\W^\prime||\le K$; we recall that $||\cdot||$ denotes the norm on $L^2(I)$. Since $\s_t$ solves 
$(ODE)_{Lag}$, this yields that $||\ddot\s_t||\le K$; by a Taylor development, this implies that
$$\norm{
\frac{\s_{t+h}-\s_t}{h}-\dot\s_t
}\le K|h|  .  \eqno (2.6)$$
The last two formulas and a Taylor development imply the first inequality below; the equality comes from an integration by parts; the last inequality comes again from (2.6).
$$\hat U(t+h,N)-\hat U(t,M)\le
h[\L_c(t,\s_t,\dot\s_t)+\a(c)]+$$
$$\int^{t+h}_0\left[
\inn{\dot\s_s-c}{\frac{(N-M)-(\s_{t+h}-\s_t)}{t+h}}
-\inn{
\V^\prime(s,\s_s)+
\W^\prime(\s_s)}{\frac{[(N-M)-(\s_{t+h}-\s_t)]s}{t+h}
}
\right]\dr s+$$
$$K(h^2+||N-M||^2)=$$
$$h[\L_c(t,\s_t,\dot\s_t)+\a(c)]+
\inn{\dot\s_{t+h}-c}{-\s_{t+h}+\s_t+(N-M)}+K(h^2+||N-M||^2)\le$$
$$h[\L_c(t,\s_t,\dot\s_t)+\a(c)]-
h\inn{\dot\s_t-c}{\dot\s_t}+\inn{\dot\s_t-c}{N-M}+
2K(h^2+||N-M||^2)  .   $$
Since 
$$-\H_0(t,\s_t,\dot\s_t)=
\L_c(t,\s_t,\dot\s_t)-\inn{\dot\s_t-c}{\dot\s_t}  ,  $$
the last formula implies that
$$\hat U(t+h,N)-\hat U(t,M)\le 
-h[\H_0(t,\s_t,\dot\s_t)-\a(c)]+\inn{\dot\s_t-c}{N-M}+
2K(h^2+||N-M||^2) .
\eqno (2.7)$$ 

To prove differentiability and (2.3), we need an inequality opposite to (2.7). We let $\s$ be as above, the minimizing curve for $\hat U(t,M)$; by hypothesis, $\s$ is unique. We note that point 4) of proposition 2.1 holds for 
$\g(0)\in Mon_3$, with the same proof. In other words, if  
$N\in Mon_3$ we can find  $\s^{h,N}$ minimal for $\hat U(t+h,N)$ such that  $\s^{h,N}_{t+h}=N$; moreover, $\s^{h,N}_s\in Mon_3$ for $0\le s\le t+h$. We set
$$\hat\s_s=\s^{h,N}_s+\frac{M-N}{t}s+
\frac{\s^{h,N}_{t+h}-\s^{h,N}_t}{t}s$$
and we see that $\hat\s_t=M$, $\hat\s_0=\s_0^{h,N}$. With the same calculations of (2.7), we get that
$$\hat U(t,M)-\hat U(t+h,N)\le$$
$$\int^t_0[\L_c(s,\hat\s_s,\dot{\hat\s}_s)+\a(c)]\dr s-
\int^{t+h}_0[\L_c(s,\s^{h,N}_s,\dot\s^{h,N}_s)+\a(c)]\dr s\le$$
$$-\inn{\dot\s^{h,N}_t-c}{N-M}+
h[\H_0(t,N,\dot\s^{h,N}_{t+h}) -\a(c)]
+K(h^2+||M-N||^2)  .  \eqno (2.8)$$
We forego the easy proof ([5]) that, if $(t+h,N)$ belongs to a ball centered in $(t,M)$, we have a uniform bound 
$$\int^t_0\norm{\dot\s^{h,N}_s}^2\dr s\le C_1 . \eqno (2.9)    $$
In particular, $\s^{h,N}_s$ is uniformly $\2$-H\"older for $|h|\le 1$ and $N\in Mon_3$, $\norm{M-N}\le 1$.

We assert that the uniform Holderianity of $\s^{h,N}_s$ implies the following: if $N\in Mon_3$ and $|h|+||M-N||<\d$, then
$$||\s^{h,N}-\s||_{C^0([0,t],Mon_3)}<\e(\d)
\txt{with} \e(\d)\tends 0
\txt{as} \d\tends 0  .  $$
It suffices to show that, if $(h_k,N_k)\tends(0,M)$ in 
$\R\times Mon_3$, then, up to subsequences, 
$\s^{h_k,N_k}\tends\s$ in $C^0([0,t],Mon_3)$. We show this fact. 

Since $\fun{\s^{h_k,N_k}}{[0,t]}{Mon_3}$ and $Mon_3$ is locally compact, we can use (2.9) and Ascoli-Arzel\`a\ as in [5] to get that, up to subsequences, $\s^{h_k,N_k}\tends\s^1$ in $C^0([0,t],Mon_3)$. Since $\s^{h_k,N_k}$ minimizes in (2.1), we easily see ([5]) that 
$\s^1$ minimizes (2.1) at $(t,M)$. By our hypotheses, $\s$ is the only minimizer; this yields that $\s^1=\s$. In other words, 
$\s^{h,N}\tends\s$ in $C^0([0,t],Mon_3)$ as $(t+h,N)\tends(t,M)$; since $\s^{h,N}$ satisfies $(ODE)_{Lag}$, it follows that 
$\s^{h,N}\tends\s$ in $C^2([0,t],Mon_3)$. This fact and (2.8) imply that, for $N\in Mon_3$,
$$\hat U(t,M)-\hat U(t+h,N)\le$$
$$-\inn{\dot\s_{t}-c}{N-M}+
h[\H_0(t,M,\dot\s_{t})-\a(c)]+
\e(\norm{M-N}+|h|)\cdot(\norm{M-N}+|h|)  $$
where $\e(\g)\tends 0$ as $\g\tends 0$.   

The last formula, together with (2.7), implies that $\hat U$ is lazily differentiable and that 
$$\partial_M\hat U(t,M)=(\dot\s_t-c),\qquad
\partial_t\hat U(t,M)=-{\cal H}_0(t,\s_t,\dot\s_t)+\a(c)  .  $$
Since $\s_t=M$, (2.3) holds.

\fin

\lem{2.6} There is $K\ge 0$ such that, for any $U\in C_\Group(\T)$, the function $T^-_1U$ is $K$-quasiconcave. In other words,  there is $K\ge 0$ such that the map 
$\Phi_K$
$$\fun{\Phi_K}{L^2(I)}{\R}, \qquad
\fun{\Phi_K}{M}{(T^-_1U)(M)-\frac{K}{2} \norm{M}^2}$$
is concave.

\proof We define a Lagrangian on $S^1\times (S^1)^n\times\R^n$ by
$$L_{n,c}(t,q,\dot q)=\frac{1}{n}\sum_{i=1}^n(
\2 |\dot q_i|^2
-c\dot q_i
)   -
\frac{1}{n}\sum_{i=1}^nV(t,q_i)-
\frac{1}{2n^2}\sum_{i,j=1}^n W(q_i-q_j)    $$
where $q=(q_1,\dots,q_n)$. 
This is the Lagrangian for the Vlasov equation with $n$ particles, each of mass $\frac{1}{n}$; its value function is
$$\hat u_n(x)=
\min \{
\int_0^1[L_{n,c}(t,q,\dot q)+\a(c)]\dt+U(D_nq(0))\st
q\in AC([0,1],\R^n), \quad q(1)=x
\}        \eqno (2.10)$$
where the operator $D_n$ has been defined at the end of section 1, and $x=(x_1,\dots,x_n)$.

Since $L_{n,c}$ is a finite-dimensional Lagrangian, the minimum above is attained by Tonelli's theorem. Let $q$ be minimal in the definition of $\hat u_n(x)$; for $h\in\R^n$, we set
$$q^{\pm h}=q(t)\pm ht . $$
Formula (2.10) implies the first inequality below.
$$\hat u_n(x+h)+\hat u_n(x-h)-2\hat u_n(x)\le$$
$$\int_{0}^1[
L_{n,c}(t,q^h,\dot q^h)+L_{n,c}(t,q^{-h},\dot q^{-h})
-2L_{n,c}(t,q,\dot q)
]\dt   =  $$
$$\int_{0}^1\{\frac{1}{n}|h|^2-\frac{1}{2n}\sum_{i=1}^n
[
V\pprime(t,q_i(t)+\theta_i^+ h_i(t-1))+
V\pprime(t,q_i(t)-\theta_i^- h_i(t-1))
]   h_i^2  -  $$
$$\frac{1}{4n^2}\sum_{i,j=1}^n[
W\pprime(q_i(t)-q_j(t)+\theta^+_{i,j}(h_i-h_j)(t-1))+
W\pprime(q_i(t)-q_j(t)-\theta^-_{i,j}(h_i-h_j)(t-1))
]  (h_i-h_j)^2     .  $$
The equality above comes from a second order Taylor development (the constants $\theta_i^\pm$ and $\theta_{i,j}^\pm$ belong to $(0,1)$ and depend on $t$); since 
$$|V^\prime{}^\prime(t,x)|\le C_1,\quad
|W^\prime{}^\prime(x)|\le C_1,
\txt{and}
(h_i-h_j)^2\le 2h_i^2+2h_j^2  ,  $$
we get that
$$\hat u_n(x+h)+\hat u_n(x-h)-2\hat u_n(x)\le
\frac{1}{n}|h|^2+
\frac{1}{2n}\sum_{i=1}^n 2C_1|h_i|^2+
\frac{1}{4n^2}\sum_{i,j=1}^n4C_1(h_i^2+h_j^2)=
\frac{K}{n}|h|^2    $$
where we have denoted by $|\cdot|$ the euclidean norm in 
$\R^n$. It is well-known that the formula above implies that the function from $\R^n$ to $\R$
$$\fun{}{x}{\hat u_n(x)-\2\frac{K}{n}|x|^2}$$
is concave. By the definition of the operator 
$\fun{D_n}{\R^n}{L^2(I)}$, we have that 
$\frac{1}{\sqrt{n}}|q|=||D_nq||$; thus, the formula above says that the function from $L^2(I)$ to $\R$
$$\fun{}{M}{\hat u_n(P_nM)-\frac{K}{2}||P_nM||^2}$$
is concave. The thesis follows from this and from the fact, proven in [5], that, if $M\in L^2(I)$, then
$$\hat u_n(P_nM)\tends (T_1^- U)(M)  \txt{as}n\tends+\infty . $$

\fin

\noindent{\bf Definition.}  Let $U\in C_{\Group}(\T)$ be $c$-dominated; we define $A_U$ as the set of the $M\in Mon$ for which there is 
$\g\in AC_{mon}([-1,1])$ with $\g_0=M$ and
$$U(\g_1)-U(\g_{-1})=
\int_{-1}^1[
\L_c(s,\g_s,\dot\g_s)+\a(c)
]\dr s   .   \eqno (2.11)$$

\thm{2.7} There is a constant $A>0$ such that, if $U$ is 
$c$-dominated and $M\in A_U$, then the following holds.

\noindent 1) $A_U$ is closed in $Mon$.

\noindent 2) If $M\in A_U$, then 
$U(M)=(T^-_1 U)(M)=(T^+_{-1} U)(M)$.

\noindent 3) If $\g$ is as in (2.11), then $\g|_{[-1,0]}$ is the unique curve on which the $\inf$ in the definition of $(T^-_1U)(M)$ is attained; analogously, $\g|_{[0,1]}$ is the unique curve on which the $\sup$ in the definition of $(T^+_{-1}U)(M)$ is attained.

\noindent 4) Let us call $\g^M$ the curve which satisfies (2.11) and $\g^M_0=M$; we recall that, by point 3), $\g^M$ is unique. Then, the map $\fun{}{M}{\dot\g^M_0}$ is continuous.

\noindent 5) For $M\in A_U$, $U$ is Fr\'echet differentiable and 
$d_MU=\dot\g^M_0-c$. Moreover, the map
$$\fun{}{A_U}{L^2(I)\times L^2(I)},\qquad
\fun{}{M}{(M,d_MU)}$$
is Lipschitz with Lipschitz inverse.

\proof We only sketch the proof, which is identical to theorem 4.5.5 of [7]. We begin with point 1). Let $M^n\in A_U$ and let 
$M^n\tends M$ in $L^2(I)$; let $\g^n\in AC_{mon}([-1,1])$ be a curve satisfying (2.11) with $\g^n_0=M^n$. We shall prove that 
$\g^n$ converges to a curve $\g$ which satisfies (2.11) and such that $\g_0=M$.

Since $\g^n$ is calibrating, it is $c$-minimal; this implies in a standard way (see lemma 3.4 below for a proof) that there is 
$C_1>0$ such that
$$\sup_{t\in(-1,1)}||\dot\g^n_t||\le C_1\qquad
\forall n\in\N  .  $$
As a consequence, $\fun{\g^n}{[-1,1]}{Mon}$ is equilipschitz and, since $\g^n_0=M^n$ is bounded, $\g^n$ is equibounded too. Since $Mon$ is a locally compact subset of $L^2(I)$, we get by 
Ascoli-Arzel\`a\ that, up to subsequences, 
$\g^n\tends\g$ in $C^0([-1,1],L^2(I))$, and that $\g_0=M$. Using the fact that $\g_n$ solves $(ODE)_{Lag}$, we see that 
$\ddot\g_n\tends\ddot\g$ in $C^0([-1,1],L^2(I))$; by the interpolation inequalities, $\g^n\tends\g$ in $C^2([-1,1],L^2(I))$. Now, the action functional is continuous under convergence in $C^1$; together with the fact that $U$ is continuous and that $\g^n$ satisfies (2.11), this implies that
$$U(\g_1)-U(\g_{-1})=
\int_{-1}^1[
\L_c(t,\g_t,\dot\g_t)+\a(c)
]\dt  ,  $$
proving point 1).

We prove point 2). Since $U$ is $c$-dominated, point 1) of lemma 2.3 implies that
$$U(M)\le (T^-_1  U)(M)\txt{and}
(T^+_{-1}  U)(M)\le U(M)
\txt{for all} M\in L^2(I) .  \eqno (2.12)$$
Let now $M\in A_U$ and let $\g$ with $\g_0=M$ satisfy (2.11). We re-write (2.11) as
$$U(\g_1)-U(M)+U(M)-U(\g_{-1})=
\int_{-1}^0[
\L_c(s,\g_s,\dot\g_s)+\a(c)
]\dr s+
\int^{1}_0[
\L_c(s,\g_s,\dot\g_s)+\a(c)
]\dr s   .   $$
Since $U$ is $c$-dominated, we also have that
$$\left\{
\eqalign{
{}&U(\g_1)-U(M)=U(\g_1)-U(\g_0)\le
\int_0^1[
\L_c(s,\g_s,\dot\g_s)+\a(c)
]\dr s\cr
{}&U(M)-U(\g_{-1})=U(\g_0)-U(\g_{-1})\le
\int^0_{-1}[
\L_c(s,\g_s,\dot\g_s)+\a(c)
]\dr s  .   }
\right.    $$
From the last two formulas, we get that
$$U(M)-U(\g_{-1})=\int^0_{-1}[
\L_c(s,\g_s,\dot\g_s)+\a(c)
]  \dr s  ,  \quad 
U(\g_1)-U(M)=\int_0^{1}[
\L_c(s,\g_s,\dot\g_s)+\a(c)
]  \dr s
 .   \eqno (2.13)$$
By the definitions of $(T^-_1 U)(M)$ and $(T^+_{-1} U)(M)$, the two formulas above imply respectively that, if $M\in A_U$,
$$(T^-_1 U)(M)\le U(M)
\txt{and}
U(M)\le (T^+_{-1}U)(M)  .  $$
This and (2.12) prove point 2).

We prove point 3). By point 2), if $M\in A_U$, then 
$(T^-_1U)(M)=U(M)$; thus, it suffices to prove that any curve 
$\tilde\g$ with $\tilde\g_0=M$ and
$$(T^-_1U)(\tilde\g_0)-U(\tilde\g_{-1})=
U(\tilde\g_0)-U(\tilde\g_{-1})=
\int_{-1}^0[
\L_c(t,\tilde\g_t,\dot{\tilde\g}_t)
]\dt      \eqno (2.14)$$
coincides with $\g$. 

Let us suppose by contradiction that $\tilde\g\not=\g$ on $[-1,0]$; we define
$$\hat\g_t=\left\{\eqalign{
\tilde\g_t &\quad t\in[-1,0]\cr
\g_t &\quad t\in[0,1]  .  
}    \right.     $$
By (2.14) and the second formula of (2.13), it follows easily that (2.11) holds for $\hat\g$. We have seen that this implies that 
$\hat\g$ is $c$-minimal on
$[-1,1]$; in particular, it satisfies $(ODE)_{Lag}$. Now $\g$ satisfies $(ODE)_{Lag}$ for the same reason; since 
$\hat\g=\g$ on $[0,1]$, we have a contradiction with the existence and uniqueness theorem. 

We prove point 4). Let $M_n\in A_U$, and let $M_n\tends M$ in 
$L^2(I)$; point 1) implies that $M\in A_U$. Let $\g^{M_n}$ satisfy (2.11) with $\g^{M_n}_0=M_n$; we see as in the proof of point 1) that the sequence 
$\g^{M_n}\in AC_{mon}(-1,1)$ has a subsequence converging to a limit $\g$ in $C^2([-1,1],L^2(I))$. As a consequence, $\g$ satisfies (2.11) and $\g_0=M$. By the uniqueness of point 3), this implies that $\g^{M_n}\tends\g^M$ in $C^2((-1,1),L^2(I))$, proving 4).

We prove point 5). Let $M\in A_U$; the inequality below is point 1) lemma 2.3; the equality, point 2) of the present theorem.
$$U(N)\le (T^-_1U)(N)\quad\forall N\in L^2(I)
\txt{and}
U(M)= (T^-_1U)(M)  .  $$
This implies the first inequality and the equality below; in the proof of proposition 2.4, we got (2.4), i. e. the second inequality below.
$$U(N)\le(T^-_1 U)(N)\le (T^-_1U)(M)+\inn{\dot\g_0^M-c}{N-M}
+K\norm{N-M}^2=$$
$$U(M)+\inn{\dot\g_0^M-c}{N-M}
+K\norm{N-M}^2  \quad\forall N\in L^2(I) . $$
Applying the same argument to $T^+_{-1}$ with time reversed, we get that
$$U(N)\ge U(M)+\inn{\dot\g_0^M-c}{N-M}
-K\norm{N-M}^2  \quad\forall N\in L^2(I) . $$
Now a general fact (proposition 4.5.3 of [7]) implies that, if the two inequalities above hold, then $U$ is Fr\'echet differentiable at any point of $A_U$, with $d_MU=\dot\g_0^M-c$. Moreover, the map 
$\fun{}{M}{d_MU}$ is Lipschitz.

\fin

\vskip 2pc
\centerline{\bf \S 3}
\centerline{\bf The minimal measures}

\vskip 1pc

\noindent {\bf Definition.} Let $M,N\in Mon$; we say that $M\simeq N$ if 
$M-N\equiv z\in\Z$. We denote by $Mon_\Z$ the space of equivalence classes; it is easy to see that $Mon_\Z$ is compact for the topology it inherits from $Mon$ (or from $L^2(I)$, which is the same.) We shall denote by $[[M]]$ the equivalence class of 
$M$ in $Mon_\Z$: we use the double brackets to avoid confusion with the equivalence class of $M$ in ${\bf S}$, which we denoted by [M]. We denote by $\Pi$ the natural projection of $Mon$ into $Mon_\Z$. In the following, we shall work mostly on $Mon$, though we shall turn to $Mon_\Z$ in all situations in which we need compactness.

We let $B_R$ be the closed ball of radius $R$ in $L^2(I)$, with the weak topology; we endow 
$$F_R\colon=S^1\times Mon_\Z\times B_R$$
with the product topology. We see that $F_R$, being the product of compact sets, is compact; moreover, it is a metric space.

\vskip 1pc

Let $\psi_s(t,M,v)$ be the flow of $(ODE)_{Lag}$; in other words, 
$$\fun{\psi_s}{\R\times L^2(I)\times L^2(I)}{\R\times L^2(I)\times L^2(I)}, \qquad
\psi_s(t,M,v)=(t+s,\g_{t+s},\dot\g_{t+s})$$
where $\g_\tau$ solves
$$\left\{\eqalign{
\ddot\g_\tau&=-\V^\prime(\tau,\g_\tau)-\W^\prime(\g_\tau)
\txt{for}\tau\in\R\cr
\g_t&=M\cr
\dot\g_t&=v  .      
}    \right.  $$

We want to restrict this flow to a compact subset 
$K_R$ of $F_R$.

\vskip 1pc

\noindent{\bf Definition.} We define the set $K_R\subset F_R$ in the following way.
Let $(t,[[M]],v)\in F_R$ and let 
$\psi_s(t,M,v)=(t+s,\g_{t+s},\dot\g_{t+s})$; if $\g_\tau\in Mon$ and 
$||\dot\g_\tau||\le R$ for all $\tau\in\R$, we say that 
$(t,[[M]],v)\in K_R$. Note that, since 
$\psi_s(t,M+k,v)=\psi_s(t,M,v)+(0,k,0)$, the condition just stated does not depend on the choice of the representative $M$.

\vskip 1pc

We are going to see below that, for $R$ large, 
$K_R\not=\emptyset$; meanwhile, we prove the following.

\lem{3.1} $K_R$ is compact in $F_R$. 

\proof Since we saw above that $F_R$ is compact, it suffices to prove that $K_R$ is closed in $F_R$. Thus, let 
$(t_n,[[M_n]],v_n)\in K_R$, and let
$(t_n,[[M_n]],v_n)\tends(t,[[M]],v)$; we must prove that 
$(t,[[M]],v)\in K_R$.

First of all, up to adding integers, we can suppose that 
$M_n\tends M$ in $Mon$.
By the definition of $K_R$, we can find a solution $\g^n$ of 
$(ODE)_{Lag}$ such that
$(\g^n_{t_n},\dot\g^n_{t_n})=(M_n,v_n)$ and $\g^n_s\in Mon$, 
$\norm{\dot\g^n_s}\le R$ for all $s\in\R$. As a consequence, 
$\fun{\g^n}{\R}{Mon}$ is $R$-Lipschitz for all $n$; since 
$(t_n,\g^n_{t_n})\tends (t,M)$, we get that $\g_n$ is locally bounded. Now, bounded sets of $Mon$ are relatively compact, and we can apply Ascoli-Arzel\`a\ and get that, up to subsequences, 
$\g^n\tends\g$ in $C^0_{loc}(\R,Mon)$. 

Clearly, $\g_s\in Mon$ for all $s$; indeed, it is the $L^2(I)$-limit of 
$\g^n_s\in Mon$.

We note that $\sup_{s\in\R}\norm{\dot\g_s}$ is l. s. c.  for the 
$C^0_{loc}(\R,Mon)$ topology; indeed, if $\norm{\psi}\le 1$, we have that
$$|\inn{\dot\g_s}{\psi}|=
\left\vert
\frac{\dr}{\dr s}\inn{\g_s}{\psi}
\right\vert\le
\sup_s\left\vert
\frac{\dr}{\dr s}\inn{\g_s}{\psi}
\right\vert
\le
\liminf_{n\tends+\infty}\sup_{\tau\in\R}
|\frac{\dr}{\dr \tau}\inn{\g^n_\tau}{\psi}|\le R . $$
The second inequality above comes from the well-known fact that, in dimension 1, the $\sup$ norm of the derivative is l. s. c. for uniform convergence; the third one comes from the fact that 
$(\tau,\g^n_\tau,\dot\g^n_\tau)\in K_R$. Since $\psi$ is arbitrary in $B_1$ and $s$ is arbitrary in $\R$, the formula above implies
$$\sup_{s\in\R}\norm{\dot\g_s}\le R  .  $$

We prove that $\g$ is an orbit with $\g_t=M$ and $\dot\g_t=v$. Since $\g^n\tends\g$ in 
$C^0_{loc}(\R,L^2(I))$ and $\g^n$ satisfies $(ODE)_{Lag}$, it follows that $\ddot\g^n\tends\ddot\g$ in $C^0_{loc}(\R,L^2(I))$. By the usual interpolation inequalities, we get that 
$\dot\g^n\tends\dot\g$ in $C^0_{loc}(\R,L^2(I))$, which implies that 
$\dot\g_t=v$. Finally, since $\g_n\tends\g$ in $C^2$, taking limits in $(ODE)_{Lag}$, we get that $\g$ solves this equation. 

\fin

\noindent {\bf Definitions.} We define $\M_1^R$ as the set of the probability measures on $K_R$, invariant by the Euler-Lagrange flow of $\L$. We note that $\M_1^R$ is not empty, if $R$ is large enough; to show this, we recall that $W^\prime(0)=0$, because $W$ is even; thus, if $q(t)$ is an orbit of the one-particle Lagrangian
$$\fun{L}{S^1\times S^1\times\R}{\R},\qquad
L(t,q,\dot q)=\2|\dot q|^2-V(t,q)  ,  $$
then $D_1q(t)$ is an orbit of $\L$; the operator $D_1$ has been defined at the end of section 1. As a consequence, if $R$ is large enough, $\M_1^R$ contains the measures induced by 
$(t,[[D_1q(t)]],\dot q(t))$, where $q$ is a periodic orbit of $L$.

We endow $\M_1^R$  with the weak$\ast$ topology; since $K_R$ is a compact metric space, we get that $\M_1^R$ is a compact metric space too.

We also define
$$\fun{I_c}{\M_1^R}{\R},\qquad
I_c(\mu)=\int_{K_R}\L_c(t,\s,v)\dr\mu(t,[[\s]],v)  .   $$

\lem{3.2} The functional $I_c$ on $\M^R_1$ is lower semicontinuous. 

\proof We note that
$$\L_c(t,\s,v)=Cin(v)-Hom(v)-P(t,\s)$$
where
$$Cin(v)=\2\norm{v}^2,\qquad
Hom(v)=\inn{c}{v},\qquad
P(t,\s)=\V(t,\s)+\W(t,\s)  .  $$
Since $V$ and $W$ are Lipschitz, and the topology on $Mon_\Z$ is the one induced by $L^2(I)$, it is immediate that
$$\fun{P}{K_R}{\R}$$
is continuous. By the definition of the weak$\ast$ topology on 
$\M^R_1$, this implies that the map
$$\fun{}{\mu}{
\int_{K_R}P(t,\s)\dr\mu(t,[[\s]],v)
}$$
is continuous. Since we have endowed $B_R$ with the weak topology, the map $\fun{Hom}{K_R}{\R}$ is continuous; as a consequence, the map
$$\fun{}{\mu}{
\int_{K_R}Hom(v)\dr\mu(t,[[\s]],v)
}  $$
is continuous too. Let us prove that the map
$$\fun{}{\mu}{
\int_{K_R}\2\norm{v}^2\dr\mu(t,[[\s]],v)
}   \eqno (3.1)$$
is l. s. c.. To do this, we let $\{ \psi_n \}_{n\ge 1}$ be a sequence dense in $B_R$ for the strong topology of $L^2(I)$, and we define
$$\fun{g_n}{L^2(I)}{\R}$$
$$g_n(v)=\sup\{
\inn{v}{\psi_i}-\2\norm{\psi_i}^2\st
i\in(1,\dots,n)
\}   .   \eqno (3.2)$$
It is a standard fact that, if $\norm{v}\le R$, then
$$\2\norm{v}^2=\sup_{\norm{\psi}\le R}
\{
\inn{v}{\psi}-\2\norm{\psi}^2
\}   .  $$
Since $\{ \psi_i \}$ is dense in $B_R$, the last formula implies that, if $v\in B_R$, then
$$g_n(v)\nearrow\2\norm{v}^2\quad\forall v\in B_R  .  $$
Formula (3.2) and Cauchy-Schwarz imply the first inequality below; since $v$ and $\psi_1$ are in $B_R$, also the second one follows.
$$g_n(v)\ge
-||v||\cdot ||\psi_1||-\2 ||\psi_1||^2\ge
-\frac{3}{2} R^2\quad\forall v\in B_R  .  $$
Since $-\frac{3R^2}{2}\in L^1(K_R,\mu)$, we can apply monotone convergence and get that
$$\int_{K_R}\2\norm{v}^2\dr\mu(t,[[\s]],v)=
\sup_{n\ge 1}\int_{K_R}g_n(v)\dr\mu(t,[[\s]],v)   .   $$
Thus, the lower semicontinuity of (3.1) follows, if we prove that each map
$$\fun{}{\mu}{\int_{K_R}g_n(v)\dr\mu(t,[[\s]],v)}$$
is continuous.  By the definition of the weak$\ast$ topology, it suffices to prove that each function $\fun{g_n}{K_R}{\R}$ is continuous. But this is true because, by (3.2), $g_n$ is the 
$\sup$ of a finite family of maps, each of which is continuous on 
$K_R$.

\fin

The next corollary follows at once from the last two lemmas.

\cor{3.3} If $c\in\R$ and $R>0$ is so large that $K_R$ is not empty, then there is $\bar\mu\in\M_1^R$ such that
$$I_c(\bar\mu)=\inf_{\mu\in\M^R_1}I_c(\mu) . $$

\rm
\vskip 1pc

We call $c$-minimal the measures which satisfy the formula above. We want to prove, following [12] and [7], that

\noindent $\bullet$ for $R$ large, the set of the $c$-minimal measures does not depend on $R$;

\noindent $\bullet$ the orbits in the support of a $c$-minimal measure are $c$-minimal.

We need a lemma.

\lem{3.4} There is a function $\fun{R}{\R}{(0,+\infty)}$, bounded on bounded sets, such that, for any 
$c$-minimal $\s\in AC_{mon}([0,1])$, we have that
$$\sup_{t\in[0,1]}\norm{\dot\s_t}\le R(c)  .  $$

\proof Since $\s_0,\s_1\in Mon$, we have $\s_0(1-)\le \s_0(0)+1$ and $\s_1(1-)\le \s_1(0)+1$; thus, we can find $z_0,z_1\in\Z$ such that  $\s_0+z_0$ and $\s_1+z_1$ have range in $[-1,1]$. We set 
$\bar A=\s_0+z_0$, $\bar B=\s_1+z_1$ and 
$$\tilde\s_t=(1-t)\bar A+t\bar B  .  $$
We denote, as usual, by $||\cdot||_{C^0}$ the $\sup$ norm; 
an easy calculation shows that
$$C_1\colon=
\norm{V}_{C^0(S^1\times S^1)}+\norm{W}_{C^0(S^1)}\ge
\norm{\V}_{C^0(S^1\times L^2)}+
\norm{\W}_{C^0(L^2)}  .  \eqno (3.3)$$
This implies the first inequality below.
$$\int_0^1\L_c(t,\tilde\s_t,\dot{\tilde\s}_t)\dt=
\int_0^1[
\2\norm{\bar B-\bar A}^2-\inn{c}{\bar B-\bar A}-
\V(t,\tilde\s_t)-\W(\tilde\s_t)
]\dt\le$$
$$\2\norm{\bar B-\bar A}^2-\inn{c}{\bar B-\bar A}+C_1\le
\frac{4}{2}+2|c|+C_1=C_2  .  $$
The last inequality above follows from the fact that $\bar A$ and 
$\bar B$ have range in $[-1,1]$. Since $\s$ is $c$-minimal, and 
$\s_i-\tilde\s_i\in L^2_\Z$ for $i=0,1$, we get that
$$\int_0^1\L_c(t,\s_t,\dot\s_t)\dt\le
\int_0^1\L_c(t,\tilde\s_t,\dot{\tilde\s}_t)\dt\le C_2  .  $$
The first inequality below follows from Cauchy-Schwarz; the second one from the fact that $\frac{1}{4}x^2-c^2\le\2 x^2-cx$; the third one, from (3.3) and the last one from the formula above.
$$\frac{1}{4}\left(
\int_0^1\norm{\dot\s_t}\dt
\right)^2
-c^2\le
\frac{1}{4}\int_0^1\norm{\dot\s_t}^2\dt-c^2\le 
\int_0^1[
\2\norm{\dot\s_t}^2-\inn{c}{\dot\s_t}
]\dt  \le$$
$$\int_0^1\L_c(t,\s_t,\dot\s_t)\dt +C_1\le
C_2+C_1  .  $$
From this it follows that
$$\int_0^1\norm{\dot\s_t}\dt\le C_4  .  \eqno (3.4)$$
We get as in (3.3) that
$$\norm{\V^\prime}_{C^0(S^1\times L^2)}+
\norm{\W^\prime}_{C^0(L^2)}\le
\norm{V^\prime}_{C^0(S^1\times S^1)}+
\norm{W^\prime}_{C^0(S^1)} .  $$
Now $\s$, being minimal, satisfies $(ODE)_{Lag}$; by the last formula, this implies that 
$$\norm{\ddot\s_t}\le C_5\quad\forall t\in[0,1]  .  $$
The last formula and (3.4) imply the thesis.

\fin

\lem{3.5} Let $R(c)>0$ be as in lemma 3.4 and let $R\ge R(c)$. Let $\mu$ minimize $I_c$ in $\M_1^R$. Then, the following three points hold.

\noindent 1) $\mu$ a. e. $(t,[[M]],v)\in K_R$ is the initial condition of a $c$-minimal orbit. Moreover, $\mu$ induces, in a natural way, a measure $\g_0$ on $K_R\cap\{ t=0 \}$ such that the following happens. If $U$ is $c$-dominated, for 
$\g_0$ a. e. $(0,[[M]],v)$, we have that the curve 
$\s_s=\pi_{mon}\circ\psi_s(0,M,v)$ is calibrating for $U$.

\noindent 2) $\mu$ is supported in $K_{R(c)}$.

\noindent 3) For the function $\a$ defined in section 2, we have that
$$-\a(c)=I_c(\mu)=\min_{\nu\in\M_1^R}I_c(\nu)  .  $$

\proof We begin to note that point 2) follows from point 1) and lemma 3.4; before proving point 1), we sketch the proof of 3) given in [7].

We have to prove only the first equality, since the fact that that 
$\mu$ is minimal in $\M_1^R$ is one of the hypotheses. Let us consider the projections
$$\fun{\pi_{mon}}{S^1\times Mon_\Z\times L^2(I)}{Mon_\Z},\qquad
\fun{\pi_{L^2}}{S^1\times Mon_\Z\times L^2(I)}{L^2(I)},$$
$$\fun{\pi_{mon\times L^2}}{S^1\times Mon_\Z\times L^2(I)}{
Mon_\Z\times L^2(I)},\qquad
\fun{\pi_{time}}{S^1\times Mon_\Z\times L^2(I)}{S^1}
   .    $$
Let us consider $(\pi_{time})_\sharp\mu$, the marginal of $\mu$ on $S^1$; since $\mu$ is invariant by $\psi_s$, it is easy to see that $(\pi_{time})_\sharp\mu$ is translation-invariant, and thus it must coincide with the Lebesgue measure on $S^1$. As a consequence, we can disintegrate $\mu$ as 
$\mu=\L^1\otimes\g_t$, where 
$\L^1$ is the Lebesgue measure on $S^1$ and $\g_t$ is a probability measure on $Mon_\Z\times B(0,R)$. Using again the fact that $\mu$ is invariant by the flow $\psi_s$, we easily see that 
$\g_t=(\pi_{mon\times L^2}\circ\psi_t(0,\cdot,\cdot))_\sharp\g_0$; as a consequence, $\g_0$ is invariant by the time-one map 
$\fun{\Psi}{(M,v)}{\pi_{mon\times L^2}\circ\psi_1(0,M,v)}$.

Let now $\fun{U}{L^2(I)}{\R}$ be a fixed point of 
$T^-_1$; we have seen in proposition 2.2 that such a function exists. By lemma 2.3, $U$ is $c$-dominated, and thus, for 
$k\in\N$, we have
$$U\circ \pi_{mon}\circ\Psi^k(M,v)-
U\circ\pi_{mon}(M,v)\le
\int_0^k[
\L_c(t,\s_t,\dot\s_t)+\a(c)
]\dt   \eqno (3.5)$$
for every $\s\in AC_{mon}([0,k])$ with 
$\s_k=\pi_{mon}\circ\Psi^k(M,v)$ and $\s_0=M$. We let 
$\s^{M,v}_t=\pi_{mon}\circ\psi_t(0,M,v)$; we consider (3.5) for 
$k\ge 1$ and $\s=\s^{M,v}$; we integrate it under $\g_0$ and we get the inequality below.
$$0=\int_{Mon_\Z\times B_R}[
U\circ\pi_{mon}\circ\Psi^k(M,v)-
U\circ\pi_{mon}(M,v)
]    \dr\g_0([[M]],v)\le$$
$$\int_{Mon_\Z\times B_R}\dr\g_0([[M]],v)
\int_0^k[
\L_c(\psi_t(0,M,v))+\a(c)
]\dt=$$
$$\int_0^k\dt
\int_{Mon_\Z\times B_R}[
\L_c(t,\tilde M,\tilde v)+\a(c)
]\dr\g_t([[\tilde M]],\tilde v)=
k\int_{K_R}[
\L_c(t,\tilde M,\tilde v)+\a(c)
]\dr\mu(t,\tilde M,\tilde v)  .  \eqno (3.6)$$
The first equality above follows because $\g_0$ is invariant by the time-one map $\Psi$, the second one because 
$\g_t=\psi_t(0,\cdot,\cdot)_\sharp\g_0$, and the third one because 
$\mu=\L^1\otimes\g_t$.

Now (3.6) implies that, for $\mu$ $c$-minimal,
$$I_c(\mu)\ge-\a(c)  .  \eqno (3.7)$$
We want to prove the opposite inequality. Let $M\in Mon$, we recall from proposition 2.2 that there is $\s\in AC_{mon}((-\infty,0])$ with 
$\s_0=M$ such that, for any $k\in\N$,
$$U(M)-U(\s_{-2k})=
\int_{-2k}^0[
\L_c(t,\s_t,\dot\s_t)+\a(c)
]\dt   .   \eqno (3.8)$$
Now we use the Krylov-Bogoljubov argument: we consider the map
$$\fun{\Phi_k}{[-k,k]}{F_R},\qquad
\fun{\Phi_k}{t}{(t\mod 1,\s_{t-k}\mod 1,\dot\s_{t-k})}$$
and the probability measure $\mu_k=(\Phi_{k})_\sharp\nu_k$, where $\nu_k$ is the Lebesgue measure on $[-k,k]$ normalized to 1. 

Since $\s$ is $c$-minimal on $(-\infty,0]$, lemma 3.4 implies that 
$\Phi_k([-k,k])\in F_{R(c)}\subset F_R$ for $k\ge 1$.

This implies that $\mu_k$ is supported in the compact set $F_R$; thus, up to subsequences, $\mu_k$ converges weak$\ast$ to a probability measure $\bar\mu$ on $F_R$. We assert that 
$\bar\mu\in{\cal M}_1^R$, i. e. that $\bar\mu$ is invariant and supported on $K_R$. The Kryolov-Bogolyubov construction implies in a standard way that $\bar\mu$ is invariant; moreover, 
$\bar\mu$ is supported on the limits of the orbits $\s_{t-k}$; but 
$\s_{t-k}\in Mon$ for $t\in(-\infty,k]$, and thus any of its limits 
$\tilde\s_t$ belongs to $Mon$ for all $t\in\R$.

This and lemma 3.2 imply the inequality below.
$$I_c(\bar\mu)\le\liminf_{k\tends+\infty}I_c(\mu_k)=
\liminf_{k\tends+\infty}\frac{1}{2k}\int_{-2k}^0\L_c(t,\s_t,\dot\s_t)\dt=
\liminf_{k\tends+\infty}\frac{1}{2k}[
U(M)-U(\s_{-2k})-2k\a(c)
]  =   -\a(c)  .  $$
The first equality above comes from the definition of $\mu_k$, 
the second one comes from (3.8) and the third one from the fact, which we saw at the beginning of section 2, that $U$ is bounded. Since $\bar\mu$ is an invariant probability measure on $K_R$, the last formula and (3.7) imply point 3).

By point 3), for $k\in\N$ formula (3.6) collapses to
$$0=\int_{K_R\cap\{ t=0 \}}[
U\circ\pi_{mon}\circ\Psi^k(M,v)-
U\circ\pi_{mon}(M,v)
]    \dr\g_0([[M]],v)=$$
$$\int_{K_R\cap\{ t=0 \}}\dr\g_0([[M]],v)
\int_0^k[
\L_c(t,\s_t^{M,v},\dot\s_t^{M,v})+\a(c)
]\dt  .   $$
This and (3.5) imply that, for all $k\in\N$ and $\g_0$ a. e. 
$([[M]],v)$, 
$$U\circ\pi_{mon}\circ\Psi^k(0,M,v)-U(M)=
\int_0^k[
\L_c(t,\s^{M,v}_t,\dot\s^{M,v}_t)+\a(c)
]\dt  .  \eqno (3.9)$$
We have seen that, since $U$ is $c$-dominated, this implies that 
$\s^{M,v}$ is $c$-minimal for a. e. $([[M]],v)$; but this is point 1).

\fin

Now we briefly define, following [12], the two "conjugate mean actions" $\a$ and $\b$.

For starters, we define the rotation number of $\mu\in\M_1^R$ in the standard way, by duality with the equivariant homology of 
$L^2(I)$. We recall from proposition 1.1 that, if $S\in C^1(L^2(I))$ and $\dr S$ is $L^2_\Z$ and $\Group$-equivariant, then
$$\dr S=c+\dr s$$
with $c\in\R$; the function $s$, which belongs to $C^1(L^2(I))$, is $L^2_\Z$ and is $\Group$-equivariant. Let $\mu\in\M_1^R$; as in [12], the ergodic theorem implies the first equality below.
$$\int_{K_R}\inn{\dr_M s}{v}\dr\mu(t,[[M]],v)=
\int_{K_R}\dr\mu(t,[[M]],v)
\lim_{n\tends+\infty}\frac{1}{n}\int_0^n
\frac{\dr}{\dr\tau} s(\pi_{mon}\circ\psi_\tau(t,M,v))\dr\tau
=$$
$$\int_{K_R}
\lim_{n\tends+\infty}\frac{1}{n}[
s\circ\pi_{mon}\circ\psi_n(t,M,v)-s\circ\pi_{mon}(t,M,v)
]\dr\mu(t,[[M]],v)=0  .  $$
The last equality above comes from the fact that any
$s\in C_{\Group}(\T)$ is bounded; we saw this right at the beginning of section 2. As a consequence,
$$\int_{K_R}\inn{c+\dr_M s}{v}\dr\mu(t,[[M]],v)$$
depends only on $c\in\R$. If we define $\r(\mu)$ as
$$\rho(\mu)=\int_{K_R}\inn{1}{v}\dr\mu(t,[[M]],v)  ,  \eqno (3.10)$$
we see by the formula above that
$$\int_{K_R}\inn{c+\dr_M s}{v}\dr\mu(t,[[M]],v)=
c\cdot\r(\mu)$$
for all $c\in\R$ and $s\in C^1(L^2(I))$, $L^2_\Z$ and 
$\Group$-equivariant. 
One can look on $\r(\mu)$ as on the "mean number of turns of all the particles around $S^1$"; indeed, by the ergodic theorem, (3.10) implies the first equality below.
$$\r(\mu)=
\int_{K_R}\dr\mu(t,[[M]],v)
\lim_{n\tends+\infty}\frac{1}{n}\int_0^n   \frac{\dr}{\dr\tau}
\inn{1}{\pi_{mon}\circ\psi_\tau(t,M,v)}\dr\tau=$$
$$\int_{K_R}\dr\mu(t,[[M]],v)
\lim_{n\tends+\infty}\frac{1}{n}
\inn{1}{[\pi_{mon}\circ\psi_n(t,M,v)-M]}  .  $$
Now, since $\fun{}{x}{\pi_{mon}\circ\psi_n(t,M,v)(x)}$ belongs to 
$Mon$, it is easy to see that
$$\lim_{n\tends+\infty}\frac{1}{n}
[\pi_{mon}\circ\psi_n(t,M,v)(x)-M(x)]$$
does not depend on $x\in I$; actually, it is equal to 
$$\lim_{n\tends+\infty}\frac{1}{n}
\inn{1}{[\pi_{mon}\circ\psi_n(t,M,v)-M]}  ,  $$
yielding that 
$$\r(\mu)=\lim_{n\tends+\infty}\frac{1}{n}
[\pi_{mon}\circ\psi_n(t,M,v)(x)-M(x)]$$
for all $x\in I$.

Let the space ${\cal C}_1$ be as in the end of section 1; it is a standard fact (see [12]) that 
$S^1\times{\cal C}_1\times{\cal C}_1$ (the phase space of a single particle), which is invariant by the Euler-Lagrange flow of 
$\L$, contains measures of any rotation number $\r\in\R$; as a consequence, if $\r\in\R$ is given and $R>0$ is large enough, 
$\M_1^R$ contains measures of rotation number $\r$.

We define
$$\b^R(\r)=\min\{
\int_{K_R}\L_0(t,M,v)\dr\mu(t,[[M]],v)\st
\mu\in\M_1^R\txt{and}\r(\mu)=\r
\}$$
and
$$-\a^R(c)=\min\{
\int_{K_R}\L_c(t,M,v)\dr\mu(t,[[M]],v)\st
\mu\in\M_1^R
\}  .  $$
The second minimum is attained by corollary 3.3; by lemma 3.2, to prove that the first minimum is attained, it suffices to prove that the set
$$\{
\mu\in\M_1^R\st\r(\mu)=\r
\} $$
is compact. Since $\M_1^R$ is compact, it suffices to prove that 
$\fun{}{\mu}{\r(\mu)}$ is continuous for the weak$\ast$ topology on $\M_1^R$; this in turn follows from the fact that the integral on the right hand side of (3.10) is a continuous function of $\mu$; we saw this in the proof of lemma 3.2, where we called it $Hom(v)$.

By point 3) of lemma 3.5, we get that, for $R\ge R(c)$, 
$\a^R(c)=\a(c)$. By point 2) of the same lemma, the $c$-minimal measures are supported in $K_{R(c)}$. By definition, $\b^R$ is decreasing in $R$; we set
$$\b(\r)=\inf_{R>0}\b^R(\r)=\lim_{R\tends+\infty}\b^R(\r)  .  $$
It is easy to see that $\a$ and $\b$ are convex; we recall the proof, which is identical to [12], that each of them is the Legendre transform of the other one. Indeed,
$$\b^\ast(c)=\sup_\r\{
\r\cdot c-\b(\r)
\}  =
\sup_{\r,R}\{
\r\cdot c-\b_R(\r)
\} =$$
$$\sup_{\r,R}\sup
\{
\r\cdot c-\int_{K_R}\L\dr\mu
\st
\r(\mu)=\r,\quad\mu\in\M_1^R
\}=$$
$$\sup_R\sup_{\mu\in\M_1^R}
\left(
-\int_{K_R}\L_c\dr\mu
\right)=\a(c)$$
where the last but one equality comes from (3.10) and the last one the fact that $\a^R(c)=\a(c)$ for $R$ large enough. The proof that $\b$ is the Legendre transform of $\a$ is analogous.

The fact that $\a$ and $\b$ are each the Legendre transform of the other, implies that both have superlinear growth. Since $\a$ is the Legendre transform of $\b$, we have that 
$\b(\r)=c\cdot\r-\a(c)$ for any $c\in\partial\b(\r)$; as a consequence, $\b(\r)$ is attained exactly on the $c$-minimal measures, for $c\in\partial\b(\r)$; since $\b$ is superlinear, 
$\partial\b(\r)$ is compact, and thus
$$R_\r\colon=\sup_{c\in\partial\b(\r)}R(c)$$
is finite. In other words, for $R\ge R_\r$, the set of measures 
$\mu$ on $K_R$ such that $\r(\mu)=\r$ and $I_c(\mu)=\b^R(\r)$ does not depend on $R$; or $\b^R(\r)=\b(\r)$ for $R$ large enough.

We define $\hat{Mat}_c$ as the closure of the union of the supports of all the $c$-minimal measures; we define
$$\tilde{Mat}_c\colon=
(\Pi, id)^{-1}\{ \hat{Mat}_c\cap \{ t=0 \}   \}
\subset Mon\times L^2(I) ,\qquad
{Mat}_c\colon=\pi_{Mon}(\tilde{Mat}_c)\subset Mon    $$
where the projection $\Pi$ was defined at the beginning of this section. In other words, $\tilde{Mat_c}$ is the set of all initial conditions $(M,v)$ in 
$Mon\times L^2(I)$ such that $(\Pi,id)\circ\psi_s(0,M,v)$ lies in the support of a $c$-minimal measure; $Mat_c$ is what we get from this set forgetting the velocity variable.

\vskip 2pc
\centerline{\bf \S 4}
\centerline{\bf The Aubry set}
\vskip 1pc

In this section, we define the Aubry set in terms of the operators $T^-_1$ and $T^+_{-1}$; we shall check that the arguments of [7] continue to work.

\lem{4.1} If $U$ is $c$-dominated, if $M\in Mat_c$ and $n\in\N$, then
$$(T^-_{n}U)(M)=(T_{-n}^+U)(M)=U(M)  .  $$

\proof Since $U$ is $c$-dominated, $U(M)\ge(T_{-n}^+U)(M)$ by point 1) of lemma 2.3. Since $M\in Mat_c$, we have that formula (3.9) holds; now (3.9) immediately implies that 
$(T_{-n}^+U)(M)\ge U(M)$, and we are done.

\fin

\thm{4.2} If $U\in C_{\Group}(\T)$ is $c$-dominated, there are a fixed point $U^-$ of $T^-_{1}$ and a fixed point $U^+$ of 
$T^+_{-1}$ which satisfy the following points.

\noindent 1) $U(M)=U^-(M)=U^+(M)$ if $M\in Mat_c$.

\noindent 2) $U^+(M)\le U(M)\le U^-(M)$ for all $M\in L^2(I)$.

\noindent 3) $U^-$ is the smallest of the fixed points of $T^-_{1}$ which are larger than $U$, and $U^+$ is the largest of the fixed points of $T_{-1}^+$ which are smaller than $U$. In other words,

\noindent $\bullet$ if $U_1^-$ is a fixed point of $T^-_{1}$ such that 
$U\le U_1^-$, then $U^-\le U_1^-$, and

\noindent $\bullet$ if $U_1^+$ is a fixed point of $T^+_{-1}$ such that 
$U\ge U_1^+$, then $U^+\ge U_1^+$.

\noindent 4) The sequences $T^-_{n} U$ and $T^+_{-n} U$ converge to $U^-$ and $U^+$ respectively, uniformly on $L^2(I)$.

\noindent 5) If $U^-$ is a fixed point of $T_{1}^-$, then there is a fixed point $U^+$ of $T^+_{-1}$ such that $U^-=U^+$ on 
$Mat_c$; moreover, $U^+\le U^-$ on $Mon$.

\proof We only sketch the proof, since it is identical to [7].

We note that
$$T^-_{n+1}U=T_n^-\circ T_1^- U\ge T_n^-U  ,  $$
where the equality comes from the semigroup property, and the inequality from point 1) of lemma 2.3 and the fact, which we saw at the beginning of section 2, that $T_n^-$ is monotone. Thus, 
$T^-_{n}U$ is an increasing sequence. Moreover, by point 2) of proposition 2.1, $T^-_{n}U$ is $L$-Lipschitz for $dist_{{\bf S}}$, for some $L>0$ independent on $n$.  Thus, $T^-_{n}U$ quotients on the compact set ${\bf S}$ as an increasing sequence of $L$-Lipschitz functions. By lemma 4.1, $T^-_{n}U=U$ on $Mat_c$; since ${\bf S}$ is compact, and 
$T^-_{n}U$ is uniformly Lipschitz, the sequence $T^-_{n}U$ is bounded in the $\sup$ norm. Thus, $T^-_{n}U$ quotients to an increasing, bounded, uniformly Lipschitz sequence of functions on ${\bf S}$; as a result, $T^-_{n}U$ converges uniformly to a 
$L$-Lipschitz function $U^-$ on ${\bf S}$. We go back to $L^2(I)$; what we just said implies that $T^-_{n}U$ converges uniformly to $U^-$ in $L^2(I)$; since 
$dist_{{\bf S}}([u],[v])=dist_{weak}(u,v)\le||u-v||$, we get that $U^-$ is $L$-Lipschitz on $L^2(I)$. Since 
$T^-_{n}U\ge U$, we get that $U^-\ge U$. Since $T^-_{n}U=U$ on $Mat_c$, we get that $U^-=U$ on $Mat_c$. Thus, $U^-$ (and 
$U^+$, with the same proof) satisfies points 1), 2) and 4).

We saw right after the definition of $\Lambda_{c,\l}$ that the map 
$\fun{}{U}{T^-_{1}U}$ is continuous for the $\sup$ norm; this implies the second equality below, while the first and last one follow by point 4).
$$T^-_{1}U^-=T^-_{1}(\lim_{n\tends+\infty}T^-_{n}U)=
\lim_{n\tends+\infty}T^-_{(n+1)}U=U^-  .  $$
This proves that $U^-$ is a fixed point of $T^-_{1}$.

We prove 3); let $U_1^-$ be as in this point. The first equality below is point 4); the inequality is the fact, which we saw before proposition 2.2, that $T^-_{n}$ is monotone: 
$T^-_{n}(V_1)\le T^-_{n}(V_2)$ if $V_1\le V_2$.
$$U^-(M)=\lim_{n\tends+\infty}T^-_{n} U(M)\le
\lim_{n\tends+\infty}T^-_{n}U_1^-(M)=U^-_1(M)  .  $$
The last equality above follows because $U_1^-$ is a fixed point of $T^-_{1}$.

We prove 5). Let $U^-$ be a fixed point of $T^-_{-1}$; by point 4), we can build $U^+$ as the limit of 
$T_{-n}^+(U^-)$ as $n\tends+\infty$; by point 1), $U^-=U^+$ on 
$Mat_c$; applying point 2) with $U=U^-$, we get that $U^+\le U^-$. 

\fin

\lem{4.3} Let ${\cal U}$ be an open neighbourhood of 
$\hat{Mat}_c$. Then, there is $t({\cal U})>0$ with the following property. If $t\ge t({\cal U})$ and $\g\in AC_{mon}([0,t])$ is 
$c$-minimal, then there is $s\in[0,t]\cap\N$ with 
$(s,\g_s,\dot\g_s)\in{\cal U}$.

\proof The proof of this lemma is identical to [7]; essentially, it follows from the fact that, as $k\tends+\infty$, the push-forward of the normalized Lebesgue measure on $[0,k]$ by the map
$\fun{}{s}{(s,\g_s,\dot\g_s)}$ accumulates on a $c$-minimal measure. We used this fact in proving point 3) of lemma 3.5.

\fin

\prop{4.4} Let $U\in C_\Group(\T)$ be $c$-dominated. Then, there is a unique couple $(U^-,U^+)$ such that $U^-$ is a fixed point of 
$T^-_{1}$, $U^+$ is a fixed point of $T^+_{-1}$ and $U^-=U^+=U$ on $Mat_c$. Moreover, $U^+\le U^-$.

\proof Existence of the couple $(U^-,U^+)$ follows from theorem 4.2. We prove uniqueness. Let $(\tilde U^-,\tilde U^+)$ be another such couple and let $M\in Mon$; since $\tilde U^-$ is a fixed point of $T^-_{1}$, by point 2) of proposition 2.2 there is 
$\s\in AC_{loc}((-\infty,0])$ such that $\s_0=M$ and, for all 
$k\in\N$, 
$$\tilde U^-(M)-\tilde U^-(\s_{-k})=
\int_{-k}^0[
\L_c(t,\s_t,\dot\s_t)+\a(c)
]\dr t  .  $$
By lemma 4.3, there is a sequence $k_j\tends+\infty$ such that 
$\s_{-k_j}\tends N\in Mat_c$. Since $\tilde U^-$ is continuous, the formula above implies that, in the formula below, the limit on the right exists and it is equal to the expression on the left.
$$\tilde U^-(M)-\tilde U^-(N)=
\lim_{j\tends+\infty}\int_{-k_j}^0[
\L_c(t,\s_t,\dot\s_t)+\a(c)
]\dr t  .  $$
Using the fact that $U^-$ is $c$-dominated, we get that
$$U^-(M)-U^-(N)\le
\lim_{j\tends+\infty}\int_{-k_j}^0[
\L_c(t,\s_t,\dot\s_t)+\a(c)
]\dr t  .  $$
Since $N\in Mat_c$, we have that $U^-(N)=U(N)=\tilde U^-(N)$; from this and the last two formulas we get that
$$U^-(M)\le\tilde U^-(M)\qquad\forall M\in Mon $$
which implies in the usual way that 
$$U^-(M)\le\tilde U^-(M)\qquad\forall M\in L^2(I)  . $$
Exchanging the r\^oles of $\tilde U^-$ and $U^-$, we get the opposite inequality; this proves the first assertion of the lemma.

The last assertion, i. e. that $U^+\le U^-$, follows, in the obvious way, from uniqueness and point 2) of theorem 4.2.

\fin

\vskip 1pc

\noindent{\bf Definition.} A pair of functions 
$U^-,U^+\in C_\Group(\T)$ is said to be conjugate if $U^-$ is a fixed point of $T^-_{1}$, $U^+$ is a fixed point of $T^+_{-1}$ and 
$U^+=U^-$ on $Mat_c$. We denote by ${\cal D}$ the set of the couples $(U^-,U^+)$ of conjugate functions. By proposition 2.2, there is a $c$-dominated function $U$; thus, by proposition 4.4, 
${\cal D}$ is not empty.

Always by proposition 4.4, if 
$(U^-,U^+)\in{\cal D}$, then $U^+\le U^-$.

We forego the easy proof that ${\cal D}$ is closed in 
$C(L^2(I),\R)\times C(L^2(I),\R)$.

\vskip 1pc

\noindent{\bf Definition.} For $(U^-,U^+)\in{\cal D}$, we set
$$\I=\{
M\in Mon\st U^-(M)=U^+(M)
\}  .  $$
Let $(U^-,U^+)\in{\cal D}$; then, by definition of conjugate couple,
$$Mat_c\subset\I  .  $$
We note that $\Pi(\I)$ is a compact set of $Mon_\Z$; indeed, we have already seen that $Mon_\Z$ is compact; since the functions 
$U^\pm$ are continuous, $\I$ is a closed set of $Mon$, implying that $\Pi(\I)$ is a closed set of $Mon_\Z$.

\thm{4.5} Let $(U^-,U^+)\in{\cal D}$ and let $M\in\I$. Then, there is a unique $c$-minimal curve $\g\in AC_{mon}(\R)$ such that 
$\g_0=M$ and, for all $m\le n\in\Z$,
$$U^\pm(\g_n)-U^\pm(\g_m)=\int_m^n[
\L_c(t,\g_t,\dot\g_t)+\a(c)
]\dt  .   \eqno (4.1) $$
In other words, $\g$ is calibrating both for $U^-$ and for $U^+$. Moreover, $U^\pm$ is Fr\'echet differentiable at $M$ and 
$$\dr_MU^+=\dr_M U^-=-c+\dot\g_0   .   \eqno (4.2)$$

\proof Let $M\in Mon$; since $U^-$ and $U^+$ are fixed points of $T^-_{-1}$ and 
$T^+_1$ respectively, we can apply point 2) of proposition 2.2 and get that there are two minimal curves, $\g^-\in AC_{loc}((-\infty,0])$ and $\g^+\in AC_{loc}([0,+\infty))$, such that
$$\g_0^-=\g_0^+=M$$
and, for any $n\in\N$ and $-m\in\N$,
$$\left\{\eqalign{
U^-(\g_0^-)-U^-(\g^-_{m})&=
\int_{m}^0[
\L_c(t,\g^-_t,\dot\g^-_t)+\a(c)
]\dt\cr
U^+(\g_n^+)-U^+(\g^+_{0})&=
\int^{n}_0[
\L_c(t,\g^+_t,\dot\g^+_t)+\a(c)
]\dt  .   
}      \right.     \eqno (4.3)    $$
We define
$$\g_t=\left\{\eqalign{
\g^-_t&\quad t\le 0\cr
\g^+_t&\quad t\ge 0  
}    \right.     $$
and we get, by (4.3), that
$$\int_{m}^n[
\L_c(t,\g_t,\dot\g_t)+\a(c)
]\dt=
U^+(\g_n)-U^-(\g_{m})+[U^-(\g_0)-U^+(\g_0)]  . \eqno (4.4) $$

We prove that, if $M\in\I$, then $\g$ satisfies (4.1); clearly, up to integer translations, we can always suppose that $m<0<n$. The first inequality below comes from the fact that 
$U^-\ge U^+$; the first equality comes from the fact that 
$\g_0=M\in\I$; the second one comes from (4.4). The last inequality comes from the fact that $U^-$ is $c$-dominated.
$$U^-(\g_n)-U^-(\g_{m})\ge
U^+(\g_n)-U^-(\g_{m})=
U^+(\g_n)-U^-(\g_{m})+[U^-(\g_0)-U^+(\g_0)]=$$
$$\int_{m}^n[
\L_c(t,\g_t,\dot\g_t)+\a(c)
]   \dt\ge
U^-(\g_n)-U^-(\g_{m})  .  $$
This formula implies (4.1) for $U^-$; the proof for $U^+$ is analogous.

We saw above that, if $U$ is $c$-dominated and $\g$ satisfies (4.1), i. e. it is calibrating, then $\g$ is $c$-minimal. This gives existence.

We prove uniqueness. Let $\tilde\g$ be any curve such that 
$\tilde\g_0=M$ and such that (4.1) holds. If we define
$$\hat\g_t=\left\{\eqalign{
\g_t&\quad t\le 0\cr
\tilde\g_t&\quad t\ge 0  
}    \right.     $$
we see as above that $\hat\g$ satisfies (4.1) and thus it is 
$c$-minimal; since $c$-minimal curves are $C^2$, we get that 
$\dot\g_0=\dot{\tilde\g}_0$; since both curves satisfy 
$(ODE)_{Lag}$, we get that $\tilde\g=\g$.

Formula (4.2) comes from (4.1) and point 5) of theorem 2.7. 

\fin

\vskip 1pc

\noindent{\bf Definition.} Let $(U^-,U^+)\in{\cal D}$; in view of theorem 4.5, we can define
$$\tilde\I=\{
(M,c+\dr_MU^-)\st M\in\I
\}  =
\{
(M,c+\dr_MU^+)\st M\in\I
\}  $$
where the derivatives are in the Fr\'echet sense.

\thm{4.6} 1) Let $(U^-,U^+)\in{\cal D}$. Then, the projection 
$$\fun{\pi_{mon}}{\tilde\I}{\I}$$
is bi-Lipschitz. 

\noindent 2) The set $\tilde\I$ is invariant by the time-one map $\Psi$ of the Euler-Lagrange flow of $\L$, and it contains the set $\tilde{Mat}_c$ defined at the end of section 3. Moreover, 
$(\Pi\times id)(\tilde\I)$ is compact in $Mon_\Z\times L^2(I)$; we recall that $\fun{\Pi}{Mon}{Mon_\Z}$ is the projection.

\noindent 3) If $(M,v)\in\tilde\I$, and if 
$\g_t=\pi_{mon}\circ\psi_t(0,M,v)$, then, for $m\le n\in\Z$,
$$U^\pm(\g_n)-U^\pm(\g_m)=
\int_m^n[
\L_c(t,\g_t,\dot\g_t)+\a(c)
]\dt  .  $$

\proof Let $(M,v)\in\tilde\I$ and let 
$\g_t=\pi_{mon}\circ\psi_t(0,M,v)$; this mean that $\dot\g_0$ satisfies formula (4.2); by the uniqueness part of theorem 4.5, it satisfies (4.1) too, and this yields point 3).

Since $\g_0=M$, setting $m=-1$ and $n=1$ in point 3) of the present theorem, we see that $M\in A_{U^-}$; the set $A_{U^-}$ has been defined before theorem 2.7. Since $M$ is arbitrary in 
$\I$, we get that $\I\subset A_{U^-}$. Point 1) follows by this and point 5) of theorem 2.7. 

Point 2): the fact that $(\Pi\times id)(\tilde\I)$ is compact follows from point 1) and the fact that $\Pi(\I)$ is compact, which we proved just before theorem 4.5. 

We prove that $\tilde\I$ is invariant by $\Psi$. Let $\g$ be as in point 3); we have that
$$\left\{\eqalign{
U^-(\g_1)-U^-(\g_{-n})&=
\int_{-n}^1[
\L_c(t,\g_t,\dot\g_t)+\a(c)
]\dt\cr
U^+(\g_n)-U^+(\g_{1})&=
\int^{n}_1[
\L_c(t,\g_t,\dot\g_t)+\a(c)
]\dt     .
}    \right.$$
Let us suppose by contradiction that $\g_1\not\in\I$, i. e. that 
$U^-(\g_1)-U^+(\g_1)>0$; summing the two formulas above, this implies that
$$U^+(\g_n)-U^-(\g_{-n})<
\int_{-n}^n[
\L_c(t,\g_t,\dot\g_t)+\a(c)
]\dt  .  $$
On the other side, since $\g_0\in\I$, we have that 
$U^-(\g_0)-U^+(\g_0)=0$; arguing as above, this implies that
$$U^+(\g_n)-U^-(\g_{-n})=
\int_{-n}^n[
\L_c(t,\g_t,\dot\g_t)+\a(c)
]\dt  .  $$
This contradiction proves that $\g_1\in\I$; since
$$-c+\dot\g_1=\dr_{\g_1}U^+=\dr_{\g_1}U^-$$
by (4.2), we get that that $\tilde\I$ is invariant by $\Psi$. 

The fact that $Mat_c\subset\I$ follows from the definition of conjugate pair; to prove that $\tilde{Mat}_c\subset\tilde\I$, we recall that, in formula (3.9), we have shown that, if 
$(M,v)\in\tilde{Mat}_c$ and $(\g_0,\dot\g_0)=(M,v)$, then $\g$ satisfies (4.1); by the uniqueness of theorem 4.5, we get that 
$\dot\g_0=c+\dr_M U^-(\g_0)$, i. e. that 
$\tilde{Mat}_c\subset\tilde\I$.

\fin

\vskip 1pc

\noindent {\bf Definition.} We define the Aubry set ${\cal A}_c$ and the Ma\~ne set ${\cal MN}_c$ in the following way.
$${\cal A}_c=\bigcap_{(U^-,U^+)\in{\cal D}}\I,\qquad
\tilde{\cal A}_c=\bigcap_{(U^-,U^+)\in{\cal D}}\tilde\I $$
$${\cal MN}_c=\bigcup_{(U^-,U^+)\in{\cal D}}\I,\qquad
\tilde{\cal MN}_c=\bigcup_{(U^-,U^+)\in{\cal D}}\tilde\I .  $$

\thm{4.7} 1) The quotiented Aubry sets $\Pi({\cal A}_c)$ and 
$(\Pi\times id)(\tilde{\cal A}_c)$ are compact; we have that $Mat_c\subset{\cal A}_c$ and
$\tilde{Mat}_c\subset\tilde{\cal A}_c$. Moreover, $\tilde{\cal A}_c$ is invariant by the time-one map $\Psi$.

\noindent 2) There is a pair $(U^-,U^+)\in{\cal D}$ such that
${\cal A}_c=\I$.

\noindent 3) The map 
$\fun{\pi_{mon}}{\tilde{\cal A}_c}{{\cal A}_c}$ is bi-Lipschitz. 

\proof By definition, each $\Pi(\I)$ is compact and contains 
$\Pi(Mat_c)$; moreover, by point 2) of theorem 4.6, each 
$(\Pi\times id)(\tilde\I)$ is compact, invariant by $\Psi$, and contains $(\Pi\times id)(\tilde{Mat}_c)$; this implies point 1).

We note that point 2) and theorem 4.6 imply point 3); actually, point 3) is also implied directly by theorem 4.6, because the restriction of a Lipschitz map to a smaller set is Lipschitz. 

We prove point 2). First of all, we restrict our conjugate couples to $Mon$, and quotient them on $Mon_\Z$; in other words, we look at them as functions in $C(Mon_\Z,\R)$. This is justified by the fact, which we saw in section 1, that any $U\in C(Mon_\Z,\R)$ can be uniquely extended to a function in $C_\Group(\T)$.

Since $Mon_\Z$ is a compact metric space, 
$C(Mon_\Z,\R)$ is separable; since ${\cal D}$ is a closed set of 
$C(Mon_\Z,\R)\times C(Mon_\Z,\R)$, we can find a dense sequence $\{ (U_n^+,U_n^-) \}_{n\ge 1}\subset{\cal D}$. Since 
$\{ U_n^\pm \}$ is a sequence of fixed points of $T^\pm_{\pm 1}$, it is equilipschitz by proposition 2.2. We note that, if 
$(U_n^+,U_n^-)\in{\cal D}$ and $a_n\in\R$, then 
$(U_n^++a_n,U_n^-+a_n)$ is a conjugate pair too; since 
$Mon_\Z$ has finite diameter, since $U_n^\pm$ is equilipschitz and $U_n^+=U_n^-$ on $Mat_c$, we can choose $a_n$ in such a way that $U_n^\pm+a_n$ is equibounded. Setting 
$\tilde U_n^\pm=U_n^\pm+a_n$, we get that the two series below converge uniformly to two Lipschitz functions on $Mon_\Z$, which we call $\tilde U^-$ and $\tilde U^+$ respectively.
$$\tilde U^-=\sum_{n\ge 1}\frac{1}{2^n}\tilde U_n^-,\qquad 
\tilde U^+=\sum_{n\ge 1}\frac{1}{2^n}\tilde U_n^+  .  \eqno (4.5)$$
Since $\tilde U_n^-=\tilde U_n^+$ on $Mat_c$, we get that 
$\tilde U^-=\tilde U^+$ on $Mat_c$.
Since $\tilde U^-$ and $\tilde U^+$ are convex combinations of 
$c$-dominated functions, it follows easily that they are $c$-dominated; by points 1) and 2) of theorem 4.2, we can find $U^-$, a fixed point of $T^-_{-1}$, satisfying $U^-\ge\tilde U^-$, and 
$U^-=\tilde U^-$ on 
$Mat_c$. Analogously, there is $U^+$, a fixed point of $T_1^+$, satisfying $U^+\le\tilde U^+$, with equality on $Mat_c$. Since 
$\tilde U^-=\tilde U^+$ on $Mat_c$, we have that $U^-=U^+$ on 
$Mat_c$, and thus $(U^-,U^+)\in{\cal D}$. As a consequence,
$${\cal A}_c\subset\I    .    \eqno (4.6)$$
On the other side, since $\{ (U_n^-,U_n^+) \}_{n\ge 1}$ is dense in 
${\cal D}$, we see that, if $M\not\in{\cal A}_c$, then
$$U_n^+(M)<U_n^-(M)$$
for at least one $n$; this implies that 
$\tilde U_n^+(M)<\tilde U_n^-(M)$ for at least one $n$. On the other side, $\tilde U_n^+\le\tilde U_n^-$ for all $n$, since 
$(\tilde U_n^-,\tilde U_n^+)$ is a conjugate pair; by (4.5), this implies that 
$$\tilde U^+(M)<\tilde U^-(M)  .  $$
We saw above that $U^+\le\tilde U^+$ and $\tilde U^-\le U^-$; thus, if $M\not\in{\cal A}_c$,
$$U^+(M)<U^-(M)  .  $$
Together with (4.6), this implies point 2).

\fin

\vskip 1pc

\noindent {\bf Definition.} Given $[[M]],[[N]]\in Mon_\Z$ and 
$n\in\N$, we define as in [13]
$$h_n([[M]],[[N]])=\min\{
\int_0^n[
\L_c(t,\g_t,\dot\g_t)+\a(c)
]\dt    
\st \s_0\in[[M]], \quad\s_n\in[[N]],
\}  $$
and
$$h_\infty([[M]],[[N]])=\liminf_{n\tends+\infty}h_n([[M]],[[N]])  .  $$
The minimum in the definition of $h_n$ is attained by an argument similar to that of point 4) of proposition 2.1. Naturally, we have to prove that $h_\infty$ is finite; for this, we refer the reader to [13], since the proof is identical.

\lem{4.8} If $(U^-,U^+)\in{\cal D}$, then
$$\forall M_-,M_+\in Mon,\quad U^-(M_-)-U^+(M_+)\le 
h_\infty([[M_-]],[[M_+]]) . $$

\proof We recall the proof of [7]. By the definition of $h_\infty$, we can find a sequence of integers $n_k\tends+\infty$ and a minimal $\g_k\in AC_{mon}([0,n_k])$ such that
$$\left\{\matrix{
& h_\infty([[M_-]],[[M_+]])=\lim_{k\tends+\infty}\int_0^{n_k}[
\L_c(t,\g^k_t,\dot\g^k_t)+\a(c)
]\dt\cr
&\g^k_0\in[[M_-]],\quad\g_{n_k}^k\in[[M_+]]  .   }
\right.    \eqno (4.7) $$
By lemma 4.3 and the fact that $\Pi(Mat_c)$ is compact in 
$Mon_\Z$, there are two integers $n_k^\prime\in[0,n_k]$ and $a_k\in\Z$ such that $\g^k_{n_k^\prime}-a_k\tends N\in Mat_c$. Since $U^-$ and $U^+$ are $c$-dominated, we have that
$$U^+(\g^k_{n_k^\prime})-U^+(M_-)\le
\int_0^{n_k^\prime}[
\L_c(t,\g^k_t,\dot\g^k_t)+\a(c)
]\dt$$
$$U^-(M_+)-U^-(\g^n_{n_k^\prime})\le
\int^{n_k}_{n_k^\prime}[
\L_c(t,\g^k_t,\dot\g^k_t)+\a(c)
]\dt  .  $$
We recall that $U^-$ and $U^+$ are $L^2_\Z$-invariant; adding the inequalities above, and letting $k\tends+\infty$, we get by (4.7) that
$$U^-(M_+)-U^-(N)+U^+(N)-U^+(M_-)\le h_\infty([[M_-]],[[M_+]])  .  $$
Since $N\in Mat_c$, the definition of ${\cal D}$ implies that 
$U^-(N)=U^+(N)$, and the thesis follows.

\fin

\thm{4.9} For $M_-,M_+\in Mon$, we have that
$$h_\infty([[M_-]],[[M_+]])=\sup_{(U^-,U^+)\in{\cal D}}
[
U^-(M_-)-U^+(M_+)
]   .  $$

\proof By lemma 4.8, we know that 
$$h_\infty([[M_-]],[[M_+]])\ge\sup_{(U^-,U^+)\in{\cal D}}
[
U^-(M_-)-U^+(M_+)
]   .   \eqno (4.8)$$
To prove the opposite inequality, we see as in theorem 5.3.6 of [7] that, for all $M_+\in Mon$, the function
$$\fun{U^-_{M_-}}{M_+}{h_\infty([[M_-]],[[M_+]])}$$
is a fixed point of $T^{-}_{1}$, while for all $M_-\in Mon$, the function 
$$\fun{U^+_{M_+}}{M_-}{h_\infty([[M_-]],[[M_+]])}$$
is a fixed point of $T_{-1}^+$. The reason for this is essentially the following: it is not hard to see that 
$\fun{Q}{M_+}{h_\infty([[M_-]],[[M_+]])}$ is $c$-dominated; moreover, the curves $\g^n$ which minimize in the definition of 

\noindent $h_n([[M_-]],[[M_+]])$ converge, up to subsequences, to a curve $\g$ calibrating for $Q$ on $(-\infty,0]$; now the assertion follows by point 2) of lemma 2.3.

Moreover, we can prove as in [7] that the conjugate function 
$U^+$ of $U^-_{M_-}$ vanishes at 
$M_-$, while the conjugate function $U^-$ of $U^+_{M_+}$ vanishes at $M_+$. Indeed, since $U^-_{M_-}$ is $c$-dominated, we can apply point 4) of theorem 4.2 and get that
$$U^+(M_-)=\lim_{n\tends+\infty}(T^+_{-n}U^-_{M_-})(M_-)=$$
$$\lim_{n\tends+\infty}\max\left\{
h_\infty([[M_-]],[[\g_n]])-\int_0^n[
\L_c(t,\g,\dot\g)+\a(c)
]\dt   \st \g_0=M_-
\right\}   .  $$
Let $\bar\g^n$ maximize in the formula above. 
For each $n$ we choose $\g^n$ minimal in the definition of 
$h_n([[M]],[[\bar\g_n^n]])$; by compactness, there is $n_k\tends+\infty$ such that $[[\g^{n_k}_{n_k}]]\tends [[N]]$. By an argument like that of point 2) of proposition 2.1, the functions $h_n$ can be shown to be $L$-Lipschitz in both variables, with the constant $L$ independent on $n$; this implies that $h_\infty$ is Lipschitz too. This, and the fact that $\g^{n_k}_{n_k}\tends N$, imply the first and third equalities below; the last one follows by our choice of 
$\g^n$.
$$\lim_{k\tends+\infty}h_\infty([[M_-]],[[\g_{n_k}^{n_k}]])=
h_\infty([[M_-]],[[N]])\le
\liminf_{k\tends+\infty}h_{n_k}([[M_-]],[[N]])=$$
$$=\liminf_{k\tends+\infty}h_{n_k}([[M_-]],[[\g^{n_k}_{n_k}]])=
\liminf_{k\tends+\infty}\int_0^{n_k}[
\L_c(t,\g^{n_k},\dot\g^{n_k})+\a(c)
\dt   .    $$
The last two formulas imply that $U^+(M_-)\le 0$. Since 
$(U_{M_-}^-,U^+)\in{\cal D}$, lemma 4.8 implies the inequality below; the equality is the definition of $U^-_{M_-}$. 
$$h_\infty([[M_-]],[[M_+]])\ge 
U^-_{M_-}(M_+)-U^+(M_-)=h_\infty([[M_-]],[[M_+]])-U^+(M_-)  .  $$
This implies that $U^+(M_-)\ge 0$, ending the proof that 
$U^+(M_-)= 0$.

Since $U^+(M_-)=0$, we get the second equality below; the first one is the definition of $U^-_{M_+}$.
$$h_\infty([[M_-]],[[M_+]])=U^-_{M_-}(M_+)=
U^-_{M_-}(M_+)-U^+(M_-)  .  $$
Since $(U^-_{M_+},U^+)\in{\cal D}$, this yields the inequality opposite to (4.8).

\fin

As an immediate consequence, we can reunite Mather's definition in [13] with Fathi's definition, which we gave before theorem 4.7.

\thm{4.10} $M\in{\cal A}_c$ iff $h_\infty(M,M)=0$.

\rm
\vskip 1pc

We forego another check, i. e. that the Ma\~ne set ${\cal MN}_c$ is the set of the $c$-minimal orbits. 

\vskip 2pc

\centerline{\bf \S 5}
\centerline{\bf Fixed points and KAM}

\vskip 1pc

Now we want to to look at the minimal orbits of $\L_c$ from another point of view, that of fixed point theory.

\vskip 1pc
\noindent {\bf Definition.} Let $\tilde\mu_{-1}$, $\tilde\mu_1$ be two Borel probability measures on $\R$, which we shall always suppose to be compactly supported. Actually, we shall only consider $\tilde\mu_{\pm 1}$ of the form 
$\tilde\mu_{\pm 1}=(\s_{\pm 1})_\sharp\nu_0$, with 
$\s_{\pm 1}\in Mon$, implying that $\tilde\mu_{\pm 1}$ is supported in an interval of length $1$.

We denote by ${\cal M}_1(\tilde\mu_1,\tilde\mu_2)$ the space of the Borel probability measures on $[-1,1]\times \R\times\R$ which satisfy the following three points.
$$\int_{[-1,1]\times \R\times\R}
(1+|v|)
   \dr\mu(t,q,v)  <  +\infty  .  \leqno i) $$

\noindent $ii$) Let $\fun{\pi}{(t,q,v)}{t}$. We ask that, if 
$\mu\in{\cal M}_1(\tilde\mu_{-1},\tilde\mu_1)$, then 
$\pi_\sharp\mu=\2\L^1$, where $\L^1$ denotes the Lebesgue measure on $[-1,1]$. In particular, $\mu$ is a probability measure, and can be disintegrated as $\mu=\2\L^1\otimes\mu_t$, with 
$\mu_t$ a measure on $\R\times\R$.

\noindent $iii$) We also ask that the elements $\mu$ of 
${\cal M}_1(\tilde\mu_{-1},\tilde\mu_1)$ are closed, i. e. for any 
$\phi\in C^1_0([-1,1]\times\R)$, we have that
$$\int_{[-1,1]\times \R\times\R}
\dr\phi(t,q)\cdot(1,v)\dr\mu(t,q,v)=
\2\int_{\R}\phi(1,q)\dr\tilde\mu_{1}(q)-
\2\int_{\R}\phi(-1,q)\dr\tilde\mu_{-1}(q)  .  \eqno (5.1) $$
In [4], the elements of 
${\cal M}_1(\tilde\mu_{-1},\tilde\mu_1)$ are called the transport measures.

\vskip 1pc

Point $i$) above essentially says that the integral on the left of (5.1) converges; point $iii$) says that $\mu$ has "boundary values" $\tilde\mu_{-1}$ at $t=-1$, and $\tilde\mu_1$ at $t=1$. As an example, consider $\s\in C^1([-1,1],L^2(I))$; if we define 
$\tilde\mu_{\pm 1}=(\s_{\pm 1})_\sharp\nu_0$, 
$\mu_t=(\s_t,\dot\s_t)_\sharp\nu_0$ and 
$\mu=\2\L^1\otimes\mu_t$, then it is easy to check that 
$\mu\in{\cal M}_1(\tilde\mu_{-1},\tilde\mu_1)$. We saw above that 
$\tilde\mu_{\pm 1}$ are supported in an interval of length $1$.

It is well-known ([4]) that we can endow 
${\cal M}_1(\tilde\mu_{-1},\tilde\mu_1)$ with a distance $d$ (called a Kantorovich-Rubinstein distance) with the following property: 
$d(\mu_n,\mu)\tends 0$ if, for any 
$\phi\in C([-1,1]\times \R\times\R)$ such that
$$\sup_{(t,q,v)}\frac{
|\phi(t,q,v)|
}{
1+|v|
}   <+\infty  ,$$
we have that
$$\int_{[-1,1]\times \R\times\R}
\phi(t,q,v)\dr\mu_n(t,q,v)\tends
\int_{[-1,1]\times \R\times\R}
\phi(t,q,v)\dr\mu(t,q,v)   .  $$
By [4], $d$ turns 
${\cal M}_1(\tilde\mu_1,\tilde\mu_2)$ into a complete metric space.

It is a standard consequence of $i$), $ii$), and $iii$) above (the proof is akin to lemma 8.1.2 of [2]) that, for any choice of the 
$C^1$ function 
$\phi$, the function
$$\fun{}{t}{
\int_{\R\times\R}\phi(q)\dr\mu_{t}(q,v)
}  $$
is absolutely continuous. In particular, the function
$$W_\mu(t,x)=\int_{\R\times\R}W(x-y)\dr\mu_t(y,v)    $$
is continuous in $t$. Since we are supposing that 
$W\in C^2(S^1)$, differentiating under the integral sign we get that 
$W_\mu\in C([-1,1],C^2(S^1))$; actually, we get that 
$|| W_\mu ||_{C([-1,1],C^2(S^1))}$ is bounded by a constant independent on $\mu$. This prompts us to define, for 
$\mu\in{\cal M}_1(\tilde\mu_{-1},\tilde\mu_1)$,
$$\fun{L_{\mu,c}}{[-1,1]\times S^1\times\R}{\R}$$
by
$$L_{\mu,c}(t,x,\dot x)=\2|\dot x|^2-c\dot x-V(t,x)-W_\mu(t,x) .  $$
An important case is that in which 
$\mu_t=(\s_t,\dot\s_t)_\sharp\nu_0$, with $\s$ $c$-minimal; we saw in section 3 that, in this case, $\s\in C^2(\R,L^2(I))$; actually, there is $C_1>0$ such that, for any $c$-minimal $\s$, 
$||\s||_{C^2(\R,L^2(I))}\le C_1$; as a consequence, 
$||W_\mu||_{C^2([-1,1]\times S^1)}\le C_2$, with $C_2$ not depending on the $c$-minimal $\s$.

To avoid proving theorems about compactness, a small haircut on transfer measures is necessary.

\vskip 1pc

\noindent{\bf Definition.} We define $A(\tilde\mu_{-1},\tilde\mu_1)$ as the smallest $R$ for which $B_R\colon=[-R,R]$ contains the supports of both $\tilde\mu_{-1}$ and $\tilde\mu_1$.

For $R\ge A(\tilde\mu_{-1},\tilde\mu_1)$, let us call 
${\cal M}_1^R(\tilde\mu_{-1},\tilde\mu_1)$ the set of the elements of ${\cal M}_1(\tilde\mu_{-1},\tilde\mu_1)$ which are supported in 
$[-1,1]\times B_R\times B_R$.  Note that 
${\cal M}_1^R(\tilde\mu_{-1},\tilde\mu_1)$ is a compact subset of 
${\cal M}_1(\tilde\mu_{-1},\tilde\mu_1)$. It follows from [3] that, for $R$ large enough, ${\cal M}_1^R(\tilde\mu_{-1},\tilde\mu_1)$ is not empty.

\lem{5.1} Let $\d\in {\cal M}_1(\tilde\mu_{-1},\tilde\mu_1)$, and let $K$ be so large that ${\cal M}_1^K(\tilde\mu_{-1},\tilde\mu_1)$ is not empty. 

\noindent 1) Then, there is 
$\bar\mu\in{\cal M}_1^K(\tilde\mu_{-1},\tilde\mu_1)$ such that
$$\int_{
[-1,1]\times \R\times \R
}
L_{\d,c}(t,q,v)\dr\bar\mu(t,q,v)=$$
$$\inf\left\{
\int_{
[-1,1]\times \R\times \R
}
L_{\d,c}(t,q,v)\dr\mu(t,q,v)
\st \mu\in
{\cal M}_1^K(\tilde\mu_{-1},\tilde\mu_1)
\right\}   .  $$

\noindent 2) The set of all the measures $\bar\mu$ which satisfy the formula above is a compact, convex set $C_\d$ of 
${\cal M}_1^K(\tilde\mu_{-1},\tilde\mu_1)$.

\noindent 3) There is $R>0$, depending on 
$A(\tilde\mu_{-1},\tilde\mu_1)$ but not on $\d$, such that, for $K\ge R$, $C_\d$ does not depend on $K$.

\proof We only sketch the standard proof of this lemma. We saw above that ${\cal M}_1^K(\tilde\mu_{-1},\tilde\mu_1)$ is compact; thus, point 1) is a standard consequence of the fact that the functional 
$$\fun{}{\mu}{
\int_{
[-1,1]\times \R\times \R
}
L_{\d,c}(t,q,v)\dr\mu(t,q,v)
}   $$
is l. s. c. (see for instance [4]). We prove point 2); since 
${\cal M}_1^K(\tilde\mu_{-1},\tilde\mu_1)$ is compact, it suffices to prove that $C_\d$ is convex and closed; this is again a consequence of the fact that the map displayed above is linear and l. s. c..

As for point 3), we recall the fact, proven in [3], that any minimal 
$\bar\mu$ is supported in a set of orbits $q$ minimal for 
$L_{\d,c}$; thus, the thesis follows if we prove that there is $R>0$, independent on $\d$, such that any minimal 
$q$, connecting a point in the support of $\tilde\mu_{-1}$ with another in the support of $\tilde\mu_1$, satisfies 
$(q(t),\dot q(t))\in B_R\times B_R$. Since $q(\pm 1)$ lie in the supports of $\tilde\mu_{\pm 1}$, i. e. in the interval 
$B_{A(\tilde\mu_{-1},\tilde\mu_1)}$, it suffices a bound on 
$\dot q$: we shall prove that $q$ satisfies $|\dot q(t)|\le C$ for a constant $C$ depending only on $A(\tilde\mu_{-1},\tilde\mu_1)$.

Actually, with the same argument of lemma 3.4, we can prove that there is $C>0$ such that, if $q$ is minimal for $L_{\d,c}$ and connects two points in $B_{A(\tilde\mu_{-1},\tilde\mu_1)}$, then 
$|\dot q|\le C$. The constant $C$ depends only on 
$|| V+W_\d ||_{C([-1,1],C^2(\T^p))}$ (which we know to be bounded independently on $\d$) and on 
$A(\tilde\mu_{-1},\tilde\mu_1)$ (the maximal distance of the points to be connected), ending the proof. 

\fin

\noindent{\bf Definition.} We settle a bit of notation: from now on, $R$ will be the constant of point 3) of the lemma above. 

If $\mu\in{\cal M}_1^R(\tilde\mu_{-1},\tilde\mu_1)$ is minimal in point 1) of lemma 5.1, we call it a minimal transfer measure for 
$L_{\d,c}$.

\vskip 1pc

Let  ${\cal C}$ denote the class of all closed, convex subsets of 
${\cal M}_1^R(\tilde\mu_{-1},\tilde\mu_1)$; by point 2) of lemma 5.1, we have a map
$$\fun{\Phi}{{\cal M}_1^R(\tilde\mu_{-1},\tilde\mu_1)}{{\cal C}} $$
which brings $\d$ into the set $C_\d$ of minimal transfer measures for $L_{\d,c}$. 

We assert that the set valued map $\Phi$ is upper semicontinuous, i. e. that, if 
$\d_n\tends\d$, if $\mu_n$ is minimal for $\L_{\d_n,c}$ and 
$\mu_n\tends\mu$, then $\mu$ is minimal for $\L_{\d,c}$. We sketch the standard proof of this; for starters, since 
${\cal M}_1^R(\tilde\mu_{-1},\tilde\mu_1)$ is closed in
${\cal M}_1(\tilde\mu_{-1},\tilde\mu_1)$, we get that 
$\mu\in {\cal M}_1^R(\tilde\mu_{-1},\tilde\mu_1)$. It is proven in [3] that the function
$$\fun{}{\mu}{
\int_{[-1,1]\times \R\times\R}[
\2|v|^2-c\cdot v-V(t,q)
]  \dr\mu(t,q,v)
}   $$
is l. s. c.. Moreover, since $\d_n\tends\d$, we have that 
$W_{\d_n}\tends W_\d$ uniformly; these two facts imply that
$$\int_{[-1,1]\times \R\times\R}
L_{\d,c}(t,q,v)\dr\mu(t,q,v)\le
\liminf_{n\tends+\infty}
\int_{[-1,1]\times \R\times\R}
L_{\d_n,c}(t,q,v)\dr\mu_n(t,q,v)  .  $$
Let us suppose by contradiction that $\mu$ is not a minimal transfer measure for $L_{\d,c}$; by the formula above, this means that there is 
$\bar\m\in{\cal M}_1^R(\tilde\m_{-1},\tilde\mu_1)$ such that
$$\int_{[-1,1]\times \R\times\R}
L_{\d,c}(t,q,v)\dr\bar\mu(t,q,v)<
\liminf_{n\tends+\infty}
\int_{[-1,1]\times \R\times\R}
L_{\d_n,c}(t,q,v)\dr\mu_n(t,q,v)  .  $$
Since $W_{\d_n}\tends W_\d$ uniformly, the formula above implies that, for $n$ large enough, the inequality below holds.
$$\int_{[-1,1]\times \R\times\R}
L_{\d_n,c}(t,q,v)\dr\bar\mu(t,q,v)=
\int_{[-1,1]\times \R\times\R}
L_{\d_,c}(t,q,v)\dr\bar\mu(t,q,v)-
\int_{[-1,1]\times \R\times\R}[
W_{\d_n}-W_\d
]\dr\bar\mu(t,q,v)<$$
$$\int_{[-1,1]\times \R\times\R}
L_{\d_n,c}(t,q,v)\dr\mu_n(t,q,v)  .  $$
This contradicts the fact that $\mu_n$ is minimal for 
$L_{\d_n,c}$, i. e. that $\mu_n\in\Phi(\d_n)$.

Since ${\cal M}_1^R(\tilde\mu_1,\tilde\mu_2)$ is compact and the map $\Phi$ is upper semicontinuous, we can apply the Ky Fan theorem ([11]) and find $\mu$ such that 
$\mu\in\Phi(\mu)$; let us gather in a set $S$ the measures $\mu$ for which $\mu\in\Phi(\mu)$. Again from the fact that $\Phi$ is u. s. c., it follows that $S$ is a closed set of 
${\cal M}_1^R(\tilde\mu_{-1},\tilde\mu_1)$; thus, it is compact, and we can find $\bar\mu\in S$ such that
$$a(\tilde\mu_1,\tilde\mu_2)\colon=
\int_{[-1,1]\times \R\times\R}
L_{\2\bar\mu,c}(t,x,v)\dr\bar\mu(t,x,v)=
\inf_{\mu\in S}
\int_{[-1,1]\times \R\times\R}
L_{\2\mu,c}(t,x,v)\dr\mu(t,x,v)   .   $$

We need a definition. 

\vskip 1pc

\noindent{\bf Definition.} We shall say that $\s$ is minimal for 
$\L_c$ if it is minimal among A. C. curves $\g$ with 
$\g_{\pm 1}=\s_{\pm 1}$. This is a weaker notion that the 
$c$-minimality of section 1, where we only required that 
$\g_{\pm 1}-\s_{\pm 1}\in L^2_\Z$. In other words, now we are considering particles on $\R$, not on $S^1$.

\vskip 1pc

The next lemma gives us the relation between the minimal transfer measures $\mu$ and the minimal paths $\s$; it can be seen as a different proof of formula (12) of [8]. Note the quirk of notation: in the definition of $a(\tilde\mu_1,\tilde\mu_2)$ we are minimizing the integral of $L_{\2\mu,c}$, but over all the minimal transfer measures $\mu$ for $L_{\mu,c}$. We shall see the reasons for this factor $\2$ in the proof below.

\lem{5.2} 1) Let $\bar\s\in AC_{mon}(-1,1)$ be minimal for $\L_c$, and let us consider the two measures on $\R$
$\tilde\mu_{-1}=(\bar\s_{-1})_\sharp\nu_0$, 
$\tilde\mu_{1}=(\bar\s_{1})_\sharp\nu_0$. Then,
$$a(\tilde\mu_{-1},\tilde\mu_1)=
\int_{-1}^1\L_c(t,\bar\s_t,\dot{\bar\s}_t)\dr t  \eqno (5.2)  $$

\noindent 2) Moreover, if $a(\tilde\mu_{-1},\tilde\mu_1)$ is attained on $\mu$, then $\mu$ is induced by a minimal parametrization 
$\s_t$; vice-versa, if $\s_t$ is a minimal parametrization, then 
$a(\tilde\mu_{-1},\tilde\mu_1)$ is attained on the measure induced by $\s_t$.

\proof  We begin with point 1). For $M_{-1},M_1\in Mon$, we define 
$$b(M_{-1},M_1)=\min\left\{
\int_{-1}^1\L_c(t,\s_t,\dot\s_t)\dt\st
\s_{-1}=M_{-1},\quad\s_1=M_1
\right\}  .  $$
Thus, we have to prove that
$$b(\s_{-1},\s_1)=a(\tilde\mu_1,\tilde\mu_2)  .  $$
We begin to show that
$$b(\bar\s_{-1},\bar\s_1)\le a(\tilde\mu_1,\tilde\mu_2)  .  
\eqno (5.3)$$
Let $\mu$ minimize in the definition of 
$a(\tilde\mu_{-1},\tilde\mu_1)$; then, 
$\mu\in S$, which implies that $\mu$ is a minimal transfer measure for $L_{\mu,c}$. 

We assert that there is a parametrization $\s\in AC_{mon}(-1,1)$ such that $\mu=\2\L^1\otimes(\s_t,\dot\s_t)_\sharp\nu_0$; note that this implies, by the definition of push-forward, that 
$$\int_{-1}^1\L_c(t,\s_t,\dot\s_t)\dt=
\int_{[-1,1]\times\R\times\R}L_{\2\mu,c}(t,x,v)\dr\mu(t,x,v)$$
from which (5.3) follows. Once we shall have proven that 
$b(\s_{-1},\s_1)=a(\tilde\mu_{-1},\tilde\mu_1)$, the formula above will yield part of point 2), i. e. that the minimal measure $\mu$ is induced by a minimal parametrization $\s$.

The proof of the assertion is essentially contained in section 4.2 of [3], which says that, if $\mu$ is a minimal transfer measure for a Lagrangian, say $L_{\mu,c}$, then $\mu$ is supported on a bunch of minimal orbits of $L_{\mu,c}$, which can be easily parametrized. 

More precisely, let $\psi_s^t$ be the Euler-Lagrange flow of 
$L_{\mu,c}$: $\psi_s^t$ brings an initial condition $(x,v)$ at time $s$ into its evolution at time $t$. By section 4.2 of [3], there is a probability measure $\tilde\mu_0$ on $\R$ and a Lipschitz function $\fun{v}{\R}{\R}$ such that, setting
$$\mu_t=(\psi_0^t)_\sharp(id,v)_\sharp\tilde\mu_0  ,  $$
then $\mu=\2\L^1\otimes\mu_t$. Take the monotone map $\s_0$ which brings $\nu_0$, the Lebesgue measure on the parameter space $[0,1]$, into $\mu_0$ and set 
$$(\s_t,\dot\s_t)=\psi_0^t\circ(id,v)\circ\s_0  .  $$
By the two formulas above, it is immediate that 
$\mu_t=(\s_t,\dot\s_t)_\sharp\nu_0$, and this proves the assertion.

We prove the inequality opposite to (5.3). Let $b(M_{-1},M_1)$ be attained on $\s$; let $\tilde\mu_{\pm 1}=(\s_{\pm 1})_\sharp\nu_0$, $\mu_t=(\s_t,\dot\s_t)_\sharp\nu_0$ and 
$\mu=\2\L^1\otimes\mu_t$. 

We begin to prove that, for $\nu_0$ a. e. $x$, the orbit 
$\fun{}{t}{\s_tx}$ is minimal, for fixed endpoints. To show this, let 
$x_0$ be a Lebesgue point for both maps 
$\fun{}{x}{(\s_{\pm 1}x,\dot\s_{\pm 1}x)}$, i. e. let 
$$\lim_{\e\tends 0+}
\frac{1}{\e}\int_{x_0-\e}^{x_0+\e}(\s_{\pm 1}x,\dot\s_{\pm 1}x)\dx=
(\s_{\pm 1}x_0,\dot\s_{\pm 1}x_0)  .  $$
We write
$$\int_{-1}^1\dt\left[
\int_I\2|\dot\s_tx|^2\dx-\int_IV(t,\s_tx)\dx-
\2\int_{I\times I}W(\s_tx-\s_ty)\dx\dr y
\right] = $$
$$\int_{-1}^1\dt\left[
\int_{I\setminus [x_0-\e,x_0+\e]}\2|\dot\s_tx|^2\dx-
\int_{I\setminus [x_0-\e,x_0+\e]}V(t,\s_tx)\dx-
\2\int_{(I\setminus [x_0-\e,x_0+\e])^2}W(\s_tx-\s_ty)\dx\dr y
\right] +\eqno (5.4)_a$$
$$\int_{-1}^1\dt\Bigg[
\int_{[x_0-\e,x_0+\e]}\2|\dot\s_tx|^2\dx-
\int_{[x_0-\e,x_0+\e]}V(t,\s_tx)\dx-
\int_{[x_0-\e,x_0+\e]\times(I\setminus [x_0-\e,x_0+\e])}
W(\s_tx-\s_ty)\dx\dr y + $$
$$\2\int_{[x_0-\e,x_0+\e]^2}
W(\s_tx-\s_ty)\dx\dr y    
\Bigg].  \eqno (5.4)_b         $$
Note that, if $x\in[x_0-\e,x_0+\e]$, the trajectory $\fun{}{t}{\s_tx}$ doesn't appear in $(5.4)_a$ of the formula above, but only in 
$(5.4)_b$; thus, the curve parametrized by 
$\s_t|_{[x_0-\e,x_0+\e]}$ minimizes the integral in $(5.4)_b$ for fixed boundary conditions. We assert that this implies that 
$\fun{}{t}{\s_tx_0}$ is minimal for $L_{\mu,c}$, endpoints fixed. 

Indeed, let us suppose by contradiction that there is 
$\fun{}{t}{q(t)}$, with the same extrema, such that
$$\int_{-1}^1L_{\mu,c}(t,q(t),\dot q(t))\dt<
\int_{-1}^1L_{c,\mu}(t,\s_tx_0,\dot\s_tx_0)\dt  .  \eqno (5.5)$$
For $\e>0$, define $\g=\g(\e)$ as the largest one among 
$\s_1(x_0+\e)-\s_1(x_0)$, $\s_1(x_0)-\s_1(x_0-\e)$, 
$\s_{-1}(x_0+\e)-\s_1(x_0)$ and $\s_{-1}(x_0)-\s_1(x_0-\e)$; since $x_0$ is a Lebesgue point of $\s_{\pm 1}$, we have that 
$\g\tends 0$ as $\e\tends 0$. 

For $x\in[x_0-\e,x_0+\e]$, we define
$$\tilde\s_tx=\left\{
\matrix{
\frac{t+1-\g}{-\g}\s_{-1}x+\frac{t+1}{\g}q(t) &-1\le t\le -1+\g\cr
q(t) & -1+\g\le t\le 1-\g\cr
\frac{t-1}{\g}q(t)+\frac{t-1+\g}{\g}\s_1x &1-\g\le t\le 1  .
}     \right.  $$
If $x\not\in[x_0-\e,x_0+\e]$, we set $\tilde\s_tx=\s_tx$. Since 
$\s_t$ solves $(ODE)_{Lag}$, for $\nu_0$ a. e. $x$ 
$\fun{}{t}{\s_tx}$ is an orbit of $L_{\mu,c}$; in particular, it is $C^2$ and depends continuously, in the $C^2$ topology, from 
the initial condition $(\s_0x,\dot\s_0x)$. Using this, the fact that $x_0$ is a Lebesgue point of $\fun{}{x}{(\s_tx,\dot\s_tx)}$ and formula (5.5), it is easy to see that
$$\int_{-1}^1\dt\Bigg[
\int_{[x_0-\e,x_0+\e]}\2|\dot{\tilde\s}_tx|^2\dx-
\int_{[x_0-\e,x_0+\e]}V(t,\tilde\s_tx)\dx-
\int_{[x_0-\e,x_0+\e]\times(I\setminus [x_0-\e,x_0+\e])}
W(\tilde\s_tx-\tilde\s_ty)\dx\dr y + $$
$$\2\int_{[x_0-\e,x_0+\e]^2}
W(\tilde\s_tx-\tilde\s_ty)\dx\dr y    
\Bigg]<$$
$$\int_{-1}^1\dt\Bigg[
\int_{[x_0-\e,x_0+\e]}\2|\dot\s_tx|^2\dx-
\int_{[x_0-\e,x_0+\e]}V(t,\s_tx)\dx-
\int_{[x_0-\e,x_0+\e]\times(I\setminus [x_0-\e,x_0+\e])}
W(\s_tx-\s_ty)\dx\dr y + $$
$$\2\int_{[x_0-\e,x_0+\e]^2}
W(\s_tx-\s_ty)\dx\dr y    
\Bigg]  $$
contradicting the fact that $\s_t|_{[x_0-\e,x_0+\e]}$ minimizes the integral in $(5.4)_b$.

Thus, for $\nu_0$ a. e. $x$ the orbit $\fun{}{t}{\s_tx}$ is minimal for 
$L_{\mu,c}$; note that we have lost the factor $\2$ in the potential of $\L_c$ (the reason for this in in formula (5.4)), and this explains the quirk of notation mentioned above.

The fact that, for $\nu_0$ a. e. $x$ the orbit $\fun{}{t}{\s_tx}$ is minimal for $L_{\mu,c}$, together with the fact that the map 
$\fun{}{\s_{-1}x}{\s_1x}$ brings $\tilde\mu_{-1}$ into $\tilde\mu_1$, implies by [3] that $\mu=\2\L^1\otimes\mu_t$ is a minimal transfer measure. Actually, this is a standard fact of transport theory: in dimension one, if minimal characteristics cannot intersect, then the unique monotone map bringing $\tilde\mu_{-1}$ into 
$\tilde\mu_1$ is a minimal transfer map.  As a consequence, we get the other half of point 2): if $\s$ is minimal, then the measure induced by $\s$ is a minimal transfer measure.

\fin

\prop{5.3} Let $\s\in AC_{mon}(\R)$ be $c$-minimal. Let us consider the set
$$A=\{
(t,\s_tx)\st t\in\R,\quad x\in I
\}\subset\R\times S^1  .  $$
Then, the map
$$\fun{\Gamma}{A}{\R},\qquad
\fun{\Gamma}{(t,\s_tx)}{\dot\s_tx}  $$
is Lipschitz.

\proof Since $\s$ is $c$-minimal, it is also minimal for $\L_c$ in the weaker sense defined above; in particular, lemma 5.2 holds.

Let $\tilde\mu_t=(\s_t)_\sharp\nu_0$, and let 
$\mu=\2\L^1\otimes(\s_t,\dot\s_t)_\sharp\nu_0$; by lemma 5.2, 
$\mu$ is a minimal transport measure for $L_{\mu,c}$ between 
$\tilde\mu_{k-1}$ and $\tilde\mu_{k+1}$, $k\in\Z$. Now we can apply the addendum in section 1.3 of of [3], which says that 
$\Gamma$ is $L$-Lipschitz for 
$t\in[k-\2,k+\2]$; the Lipschitz constant $L$ depends only on the 
$C^2$ norm of $V$ and $W_\mu$, and on the radius of the smallest ball containing the supports of $\tilde\mu_{k-1}$ and 
$\tilde\mu_{k+1}$. Since we are free to translate by an integer, $L$ depends on the diameter of the union of the supports of 
$\tilde\mu_{k-1}$ and $\tilde\mu_{k+1}$.  We note that $V$ is fixed, while $||W_\mu||_{C^2(\R\times\R)}$ is bounded, since 
$W\in C^2(S^1)$ and the $C^2$ norm of $\s_t$, which solves 
$(ODE)_{Lag}$, is bounded (we saw in lemma 3.4 that 
$\sup_{t\in\R}\norm{\dot\s_t}$ is bounded). Since $\s_t$ belongs to $Mon$ and has bounded speed in $L^2(I)$, the diameter of the union of the supports of $\tilde\mu_{k-1}$ and $\tilde\mu_{k+1}$ is bounded, uniformly in $k$. Thus, the Lipschitz constant $L$ of $\Gamma$ on $[k-\2,k+\2]$ does not depend on $k$, and the thesis follows.

\fin

We want to study the regularity of periodic minimal measures of irrational rotation number. We consider the Lagrangian
$$\L^\e(t,\s,v)=
\2\norm{v}^2-\e\V(t,\s)-\e\W(\s)$$
where $\V$ and $\W$ are defined as in section 1; the only difference is that we ask that the potentials $V$ and $W$ are 
$C^k$ for some large $k$ which we shall determine in the following.

We want to study
$$\min\{
\int_0^1
\L^\e(t,\s_t,\dot\s_t)\dt\st\s\in{\cal P}_\o
\}   \eqno (5.6)$$
where by ${\cal P}_\o$ we denote the set of those 
$\s\in AC_{mon}([0,1])$ such that 

\noindent $\bullet$ $\s_0\simeq\s_1$ in the sense of section 1, or $\s$ projects to a periodic curve on ${\bf S}$; in other words,
$$(\s_{0})_\sharp\nu_0=(\s_{1})_\sharp\nu_0   . 
\eqno (5.7)$$

\noindent $\bullet$ Moreover, we ask that the rotation number of 
$\s$ is $\o$; in other words,
$$\int_0^1\inn{1}{\dot\s_t}\dt=\o  .  \eqno (5.8)$$
We note that ${\cal P}_\o$ is closed in 
$AC_{mon}([0,1])$: we forego the easy proof that (5.7) and (5.8) are closed under uniform convergence. Moreover, ${\cal P}_\o$ is not empty, since it is easy to see that $\s_tx=x+\o t$ is periodic in ${\bf S}$ (actually, it is constant) and has rotation number $\o$. As a consequence, we were justified in writing $\min$ in (5.6).

Let $\g,\tau>0$; we say that $\o\in\R$ is $(\g,\tau)$-diophantine if
$$|\o q-p|\ge\frac{\g}{q^\tau}\txt{if}
(q,p)\in(\N\setminus \{ 0 \})\times\Z  .  $$
We want to prove the following.

\thm{5.4} Let $\o$ be $(\g,\tau)$-diophantine, and let $\s$ minimize in (5.6). Then, there is $k_0(\g,\tau)>0$ such that, if 
$V$ and $W$ are $C^k$ with $k\ge k_0(\g,\tau)$ and $\e$ is small enough, the measure on $[0,1]\times S^1\times\R$ given by
$\mu\colon=\L^1 \otimes(\s_{t},\dot\s_t)_\sharp\nu_0$ is the 
push-forward of the Lebesgue measure on $S^1\times S^1$ by a $C^1$ map $\fun{\Phi}{(t,x)}{(t,\phi_1(t,x),\phi_2(t,x))}$. 

\proof Let $\s$ be minimal in (5.6), and let 
$\mu=\L^1\otimes(\s_t,\dot\s_t)_\sharp\nu_0$; let the potential $W_\mu(x,t)$ and the Lagrangian $L_{\mu,0}(t,x,\dot x)$ be defined at the beginning of this section. Since $\s_0\simeq\s_1$ by (5.7), the definition of $W_\mu$ implies that 
$W_\mu(1,x)=W_\mu(0,x)$; thus $L_{\mu,0}$ is 1-periodic in time. 

Since $\s$ is minimal, it is a periodic solution of $(ODE)_{Lag}$; since the potentials $V$ and $W$ are $C^k$, we get that 
$\s\in C^{k+1}(S^1,L^2(I))$; as a consequence, 
$W_\mu\in C^k(S^1\times S^1)$, while 
$V\in C^k(S^1\times S^1)$ by hypothesis. Using again the fact that $\s$ is minimal, we see as in lemma 3.4 that 
$||\dot\s||_{C^0(\R,L^2(I))}$ is bounded, independently on 
$\e\in[0,1]$; differentiating in $(ODE)_{Lag}$, we get that the higher derivatives are bounded too.  Thus, $||\s||_{C^k(\R,L^2(I))}$ is bounded independently on $\e\in[0,1]$; as a consequence, 
$||W_\mu||_{C^k(\R\times S^1)}$ is bounded independently on 
$\e$. In particular, $||\e V+\e W_\mu||_{C^k(S^1\times S^1)}$ tends to zero as $\e\tends 0$; thus, by [15], for $\e$ small and $k$ large enough, $L_{\mu,0}$ has a KAM torus of rotation number 
$\o$.

We are supposing that $\s$ is minimal in (5.6); by periodicity, this implies that $\s$ is minimal, with fixed boundary conditions, on each interval $[t_0,t_0+1]$. From lemma 5.2 we gather that, for 
a. e. $x\in I$, $\s_tx$ is minimal for $L_{\mu,0}$ on each interval 
$[t_0,t_0+1]$ for fixed boundary conditions; in particular, 
$\fun{}{t}{\s_tx}$ is an orbit of $L_{\mu,0}$. This immediately implies that $\mu$ is invariant by the Euler-Lagrange flow of 
$L_{\mu,0}$. Moreover, (5.8), and the fact that $\s_t\in Mon$, imply as in section 3 that
$$\lim_{t\tends+\infty}\frac{\s_tx-\s_0x}{t}=\o
\qquad \forall x\in I . \eqno (5.9)$$

We saw above that there is $k_0(\tau,\g)$ such that, if 
$k\ge k_0(\tau,\g)$ and $\e$ is small enough, then $L_{\mu,0}$ has a KAM torus of frequency $\o$. In other words, there is a $C^1$ map $\fun{\Phi}{S^1\times S^1}{S^1\times S^1\times\R}$ such that, denoting as usual by $\psi_s$ the Euler-Lagrange flow of $L_{\mu,0}$,
$$\psi_s\circ\Phi(t,x)=\Phi(t+s,x+\o s)  ,  $$
or $\Phi$ conjugates the rotation on $S^1\times S^1$ given by 
$\fun{}{(t,x)}{(t+s,x+\o s)}$ to the Euler-Lagrange flow on the image of $\Phi$. We have to show that $\mu$ is the push-forward by $\Phi$ of the Lebesgue measure on 
$S^1\times S^1$. Since the KAM torus is conjugate to an irrational rotation, it supports just one invariant measure, i. e. the push-forward of Lebesgue. Thus, it suffices to prove that $\mu$, which we proved to be invariant, is supported on the KAM torus; equivalently, that, for each $x\in I$, the orbit $(t,\s_tx,\dot\s_tx)$ lies on the KAM torus. This is a consequence of (5.9) and of the fact that 
$\fun{}{t}{\s_t x}$ is an orbit. We explain why.

By [12] and [13], we know that the KAM torus is a graph; in other words, there is a Lipschitz map
$$\fun{v}{S^1\times S^1}{\R} $$
such that the image of $\Phi$ coincides with the graph of $v$. Moreover, the two sets
$$A_-=\{
(t,q,\dot q)\st \dot q<v(t,q)
\}  , \qquad
A_+=\{
(t,q,\dot q)\st \dot q>v(t,q)
\}  $$
are invariant by the flow $\phi_s$. 

Let us call ${\cal T}$ the KAM torus of frequency $\o$; it is standard that both $A_-$ and $A_+$ contain sequences ${\cal T}^n_-$ and 
${\cal T}^n_+$ respectively of KAM tori which, as $n\tends+\infty$, converge to ${\cal T}$. Since no orbit can cross a KAM torus, we get that, for any $z\in I$, the closure $C$ of 
$$\{
(t,\s_tz,\dot\s_tz)
\}_{t\in\R}  $$
either is contained in ${\cal T}$, or in one ot the two invariant sets  $A_\pm$, and at a finite distance from ${\cal T}$. Let us suppose by contradiction that, for some $z\in I$, $C$ is not contained in 
${\cal T}$; to fix ideas, let $C\subset A_+$.

Let us denote by $q_{x,t}(s)$ the orbit on the KAM torus such that $q_{x,t}(t)=x$. We assert two facts:

\noindent 1) if $x^\prime>x$, then there is a positive number 
$\d(x^\prime-x)$, only depending on $x^\prime-x$, such that
$$q_{x^\prime,t}(s)\ge q_{x,t}(s)+\d(x^\prime-x)
\qquad\forall s\in\R . $$

\noindent 2) There is $\e>0$, independent on $t\in\R$, such that, if $\s_tz=x$, and if 
$q_{x,t}(s)=q_{x,t}(t)+1=x+1$, then $\s_sz\ge q_{x,t}(s)+1+\e$.

Before proving 1) and 2), we show how they imply the thesis. Let $z\in I$ and $C$ be as above; we set $\s_0z=x$; by 2), we see that, if 
$q_{x,0}(s)=x+1$, then $\s_sz\ge x+1+\e$. Now we set 
$x^\prime=\s_sz$; applying again point 2), we have that, if 
$q_{x^\prime,s}(s_1)=x^\prime+1$, then 
$\s_{s+s_1}z\ge x^\prime+1+\e$. By point 1), this means that 
$\s_{s+s_1}z\ge q_{x,0}(s+s_1)+\d(\e)+\e$. Iterating, we have that
$$\s_{s+s_1+\dots+s_n}z\ge
q_{x,0}(s+s_1+\dots+s_n)+(n-1)\d(\e)+\e  .  $$
This fact implies the inequality below; the equality comes from the fact that the KAM torus has rotation number $\o$.
$$\lim_{n\tends+\infty}
\frac{\s_{s+s_1+\dots+s_n}z}{s+s_1+\dots+s_n}\ge
\lim_{n\tends+\infty}
\frac{q_{x,0}(s+s_1+\dots+s_n)}{s+s_1+\dots+s_n}+\d(\e)
=\o+\d(\e)   .  $$
We have reached a contradiction with (5.9).

We prove the two assertions above. To prove point 1), we begin to note that 
$$\Phi(t,x)=(t,\Phi_x(t,x),\Phi_v(t,x))  .  $$
Now point 1) is true for the rotation $\fun{}{(t,x)}{(t+s,x+\o s)}$, with 
$\d(x^\prime-x)=x^\prime-x$; since $\Phi$ is a conjugation, we have that $q_{x,t}(s)=\Phi_x(t+s,x+\o s)$; thus, it suffices to show that the map $\fun{}{x}{\Phi_x(t,x)}$ is strictly monotone for all $t$. This follows since, by the KAM theorem, the map 
$\fun{}{(t,x)}{(t,\Phi_x(t,x))}$ is close to the identity.

Since $C\subset A_+$ is at finite distance from ${\cal T}$, we get from the definition of $A_+$ that there is $a>0$ such that
$$\dot\s_tz\ge v(\s_tz)+a\qquad
\forall x\in I,\quad\forall t\in\R  .  \eqno (5.10)$$
Let now $x\in I$, $t\in\R$ and let $q_{x,t}(s)$ be the orbit of the KAM torus with initial conditions $q_{x,t}(t)=\s_tz=x$ and 
$\dot q_{x,t}(t)=v(t,x)$. By (5.10), $\s_sz>q_{x,t}(s)$ for $s-t$ positive and small; let $(t,T)$ be the largest interval on which 
$\s_sz>q_{x,t}(s)$. We assert that $T=+\infty$. Indeed, if $T$ were finite, then we would have 
$\s_Tz=q_{x,t}(T)$; together with $\s_sz>q_{x,t}(s)$ for 
$s\in(0,T)$, this implies that $\dot\s_Tz\le\dot q_{x,t}(T)$, contradicting (5.10). As a consequence, if $s$ is such that 
$q_{x,t}(s)=q_{x,t}(t)+1=\s_tz+1$, then $\s_sz\ge \s_tz+1+\e$, for some $\e>0$ independent on  $t\in\R$. 

\fin

\vskip 2pc
\centerline{\bf Bibliography}


\noindent [1] R. A. Adams, Sobolev spaces, Academic press, 1975, New York.

\noindent [2] L. Ambrosio, N. Gigli, G. Savar\'e, Gradient flows, Birkh\"auser, Basel, 2005.


\noindent [3] P. Bernard, B. Buffoni, Optimal mass transportation and Aubry-Mather theory, J. Eur. Math. Soc., {\bf 9}, 85-121, 2007.

\noindent [4] P. Bernard, Young measures, superposition and transport, Indiana Univ. Math. J., {\bf 57}, 247-275, 2008.

\noindent [5] U. Bessi, Chaotic orbits for a version of the Vlasov equation, Siam J. Math. Anal., {\bf 44}, 2496-2525, 2012.

\noindent [6] E. Carlen, Lectures on optimal mass transportation and certain of its applications, mimeographed notes, 2009.

\noindent [7] A. Fathi, Weak KAM theorem in Lagrangian dynamics, Fourth preliminary version, mimeographed notes, Lyon, 2003.

\noindent [8] W. Gangbo, A. Tudorascu, Lagrangian dynamics on an infinite-dimensional torus; a weak KAM theorem, Advances in Mathematics, {\bf 224}, 260-292, 2010.

\noindent [9] W. Gangbo, Hwa Kil Kim, T. Pacini, Differential forms on Wasserstein space and infinite-dimensional Hamiltonian Systems, Providence, 2010.

\noindent [10] W. Gangbo, T. Nguyen, A. Tudorascu, Hamilton-Jacobi equations in the Wasserstein space, Methods Appl. Anal., 
{\bf 15}, 155-183, 2008.

\noindent [11] Ky Fan, Fixed point and minimax theorem in locally convex topological linear spaces, Proceedings of the Nat. Acad. Sciences USA, {\bf 38}, 121-126, 1952.


\noindent [12] J. N. Mather, Action minimizing invariant measures for 
positive-definite Lagrangian Systems, Math. Zeit., {\bf 207}, 
169-207, 1991.


\noindent [13] J. N. Mather, Variational construction of connecting orbits, Ann. Inst. Fourier, {\bf 43}, 1349-1386, 1993.

\noindent [14] L. Nurbekian, Weak KAM theory on the $d$-infinite dimensional torus, Ph. D. thesis, 2012.

\noindent [15] D. Salamon, E. Zehnder, KAM theory in configuration space, Comment. Math. Helvetici, {\bf 64}, 84-132, 1989.


\end

For starters,  we suppose that $\tilde\mu_{-1}$ and 
$\tilde\mu_1$ are the sum of finitely many Dirac deltas:
$$\frac{1}{n}\sum_{i=1}^n\d_{x_i^{-}}=\tilde\mu_{-1},
\qquad
\frac{1}{n}\sum_{i=1}^n\d_{x_i^{+}}=\tilde\mu_{1}  .  $$
We can suppose that 
$$x_1^-<x_2^-<\dots<x_n^-,
\txt{and}
x_1^+<x_2^+<\dots<x_n^+   .   $$
Let $\mu$ minimize in $a(\tilde\mu_{-1},\tilde\mu_1)$; then, 
$\mu\in S$, which implies that $\mu$ is a minimal transfer measure for $L_{\mu,c}$; by [3], it is supported on a set of orbits, minimal for $L_{\mu,c}$, which don't intersect among themselves. Since we are on $\R$, this and (5.1) easily imply that 
$$\mu=\2\L^1\otimes\frac{1}{n}
\sum_{i=1}^n\d_{q_i(t)}\otimes\d_{\dot q_i(t)}  $$
with $q_i$ minimal and connecting $x_i^-$ with $x_i^+$. This fact, whose proof in this case is an exercise, is one of the basic facts of transport theory: if minimal characteristics cannot cross, then the only minimal connection between two $n$-uples of deltas, say 
$x_1^-<x_2^-<\dots<x_n^-$ and $x_1^+<x_2^+<\dots<x_n^+$, is the map which brings $x_i^-$ into $x_i^+$.

If we define 
$\s_tx=D_n(q_1(t),\dots,q_n(t))$ for the operator $D_n$ of section 1, an easy calculation implies
$$\2\int_{-1}^1\L_c(t,\s_t,\dot\s_t)\dt=
\int_{[-1,1]\times\R\times\R}L_{\2\mu,c}(t,x,\dot x)
\dr\mu(t,x,\dot x)  $$
proving (5.3) when $\tilde\mu_{\pm 1}$ has the form above. It is standard (see [4]) that $a$ is continuous in 
$(\tilde\mu_{-1},\tilde\mu_1)$, while it is proven in [5] that 
$b(\s_{-1},\s_1)$ is continuous; this implies (5.3) in the general case.

We prove the inequality opposite to (5.3); again, we start assuming that $\tilde\mu_{\pm 1}$ have the form above; in other words, for the operator $D_n$ defined in section 1,  
$\bar\s_{\pm 1}=D_n(x_1^\pm,\dots,x_n^\pm)$.
An easy computation (see again [5]) implies that, if 
$\g_t=D_n(q_1(t),\dots,q_n(t))$ is minimal with boundary conditions $P_n\bar\s_{\pm 1}$, then
$$\int_{-1}^1\L_c(t,\g_t,\dot\g_t)\dt=
\int_{-1}^1[
\frac{1}{n}\sum_{i=1}^n(
\2|\dot q_i|^2-c\dot q_i-V(t,q_i)
)-\frac{1}{2n^2}\sum_{i,j=1}^nW(q_i(t)-q_j(t))
]\dt   .  $$
Since $W$ is even, we can write the equality above as 
$$\int_{-1}^1\L_c(t,\g_t,\dot\g_t)\dt=
\int_{-1}^1[
\frac{1}{n}\sum_{i\not=j}(
\2|\dot q_i|^2-c\dot q_i-V(t,q_i)
]\dt-
\int_{-1}^1\frac{1}{2n^2}\sum_{i,l\not=j}^nW(q_i(t)-q_l(t)) \dt+$$
$$\frac{1}{n}\int_{-1}^1[
\2|\dot q_j|^2-c\dot q_j-V(t,q_j)-
\frac{1}{n}\sum_{i=1}^nW(q_i(t)-q_j(t))
]\dt  .  \eqno (5.4)$$
In the formula above, the first two terms on the right do not contain $q_j$; as a consequence, each $q_j$ must minimize the Lagrangian
$$L^n_c(t,q,\dot q)=
\2|\dot q|^2-c\dot q-V(t,q)-
\frac{1}{n}\sum_{i=1}^nW(q_i(t)-q)  .  $$
In other words, each $q_j$ is minimal for the same Lagrangian; 
now it is standard (see for instance [6]) that, since each $q_i$ connects $x_i^-$ with $x_i^+$ and is minimal for $L^n_c$, then
$$\mu_n\colon=\2\L^1\otimes\frac{1}{n}
\sum_{i=1}^n(\d_{q_i(t)}\otimes\d_{\dot q_i(t)})  \eqno (5.5)$$
is a minimal transfer measure for $L^n_c$. In the language of transport, when we connect an n-uple of deltas with another n-uple of deltas, there is not just a minimal transfer plan, but a minimal transfer map; when we are in $\R$, this is the unique monotone map bringing one measure into the other.

Now, $\mu_n$ is a minimal transfer measure for $L^n_c$, a Lagrangian which coincides with $L_{\mu_n,c}$; thus, 
$\mu_n\in\Phi(\mu_n)$; a calculation like the one we did above shows that
$$\int_{-1}^1L_{\2\mu_n,c}(t,q,\dot q)\dr\mu_n(t,q,\dot q)=
\int_{-1}^1\L_c(t,\g_t,\dot\g_t)\dt$$
and this proves the inequality opposite to (5.3) when 
$\tilde\mu_{-1}$ and 
$\tilde\mu_1$ are the sum of finitely many deltas; taking limits, we get the inequality opposite to (5.3) in the general case, proving (5.2).

Note that, in (5.4), we saw that each $q_i$ is minimal for 
$L^n_c=L_{\mu_n,c}$: we have lost the factor $\2$ in the potential of $\L_c$, and this explains the quirk of notation mentioned above. 

We prove point 2). Let $\bar\s$ be minimal on $(-1-\e,1+\e)$. We recall that in [5] it is proven that, defining $P_n$ and ${\cal C}_n$ as in section 1, 
$$\int_{-1}^1\L_c(t,\bar\s_t,\dot{\bar\s}_t)\dt=
\lim_{n\tends+\infty}\min\Bigg\{
\int_{-1}^1\L_c(t,\g_t,\dot\g_t)\dt\st$$
$$\g_{\pm 1}=P_n\s_{\pm 1},\quad
\g_t\in{\cal C}_n\quad\forall t\in[-1,1]
\Bigg\}   .    $$
The minimum on the right of the formula above is a minimum on curves arriving in ${\cal C}_n$, which is a finite-dimensional space; thus, the minimum exists by Tonelli's theorem. Let $\g^n$ be a curve minimal on the right of the formula above. By [5], up to subsequences, $\g^n$ converges uniformly to $\g$, a curve which minimizes 
$$\int_{-1}^1\L_c(t,\g_t,\dot\g_t)\dt$$
with boundary conditions $\g_{\pm 1}=\bar\s_{\pm 1}$. We recall the standard argument proving that $\g_t=\bar\s_t$ for 
$t\in[-1,1]$. Indeed, if they were different, then we could define
$$\tilde\s_t=\left\{
\eqalign{
\g_t &\qquad t\in[-1,1]\cr
\bar\s_t &\qquad t\in[-1-\e,1+\e]\setminus[-1,1]  .  
}     \right.     $$
It is easy to check that $\tilde\s$ is minimal for $\L_c$ on 
$[-1-\e,1+\e]$, endpoints fixed; however, by the existence and uniqueness theorem, it doesn't satisfy $(ODE)_{Lag}$ at 
$t=\pm 1$. By this contradiction, $\g^n$ converges uniformly to the minimal $\bar\s$.

By [5], we have that $\g^n_t=D_n(q_1(t),\dots,q_n(t))$, with 
$(q_1,\dots,q_n)$ minimal for the Lagrangian 
$$\L_c(t,D_nq,D_n\dot q)=
\frac{1}{n}\sum_i[
\2|\dot q|^2-c\dot q_i-V(t,q_i)
]  -
\frac{1}{2n^2}\sum_{i,j}W(q_i-q_j)  .  $$
Let $\mu_n$ be defined as in formula (5.5); we saw above that the minimality of $\g^n_t$ implies that $\mu_n$ is a minimal transfer measure for $L_{\mu_n,c}$; equivalently, if 
$\bar\g^n$ is ${\cal C}_n$-valued and has the same endpoints as 
$\g_n$,
$$\int_{-1}^1\dt\int_I L_{\mu_n,c}(t,\g_t^nx,\dot\g_t^nx)\dx  \le
\int_{-1}^1\dt\int_I L_{\mu_n,c}(t,\bar\g_t^nx,\dot{\bar\g}_t^nx)\dx   .   $$
Since the potential in $L_{\mu_n,c}$ converges uniformly to the potential of $L_{\bar\mu,c}$, and $\g^n\tends\bar\s$ uniformly, we get that 
$\bar\mu$ is a minimal transfer measure for $L_{\bar\mu,c}$, and that $\bar\s$ minimizes
$$\int_{-1}^1\int_I L_{\bar\mu,c}(t,\s_tx,\dot\s_tx)\dx$$
among all A. C. functions with $\s_{\pm 1}=\bar\s_{\pm 1}$, ending the proof of point 2).

We define the probability measure $\tilde\mu$ on 
$Mon_\Z\times B_R$ by
$$\tilde\mu(A)=\mu
\{
[[\psi_{t}(0,M,v)]]\st ([[M]],v)\in A, \quad t\in (-1,0]
\}  .   $$
This is actually the push-forward of $\mu$ under the map which brings $(t,\g_t,\dot\g_t)$ to $(0,\g_0,\dot\g_0)$. 
Since $\mu$ is invariant by the Euler-Lagrange flow $\psi_t$, we get that $\tilde\mu$ is invariant by the time-one map 
$\fun{\Psi}{(M,v)}{\pi_{mon\times L^2}\circ\psi_1(0,M,v)}$. Let 
$\fun{U}{L^2(I)}{\R}$ be a fixed point of 
$T^-_1$; we have seen in proposition 2.2 that such functions exist. By lemma 2.3, $U$ is $c$-dominated, and thus, for $k\in\N$, we have
$$U\circ \pi_{mon}\circ\Psi^k(M,v)-
U\circ\pi_{mon}(M,v)\le
\int_0^k[
\L_c(t,\s_t,\dot\s_t)+\a(c)
]\dt   \eqno (3.5)$$
for every $\s\in AC_{mon}([0,k])$ with 
$\s_k=\pi_{mon}\circ\Psi^k(M,v)$ and $\s_0=M$. We let 
$\s^{M,v}_t=\pi_{mon}\circ\psi_t(0,M,v)$; we consider (3.5) for 
$k\ge 1$ and $\s=\s^{M,v}$; we integrate it under $\tilde\mu$ and we get the inequality below.
$$0=\int_{K_R\cap\{ t=0 \}}[
U\circ\pi_{mon}\circ\Psi^k(M,v)-
U\circ\pi_{mon}(M,v)
]    \dr\tilde\mu([[M]],v)\le$$
$$\int_{K_R\cap\{ t=0 \}}\dr\tilde\mu([[M]],v)
\int_0^k[
\L_c(t,\s_t^{M,v},\dot\s_t^{M,v})+\a(c)
]\dt=$$
$$\int_{K_R\cap\{ t=0 \}}\dr\tilde\mu([[M]],v)
\sum_{j=0}^{k-1}\int_0^1\L_c\circ\psi_s\circ
(j,\Psi^j(M,v))\dr s=$$
$$k\int_{K_R\cap\{ t=0 \}}\dr\tilde\mu([[M]],v)\int_0^1
\L_c\circ\psi_s(0,M,v)\dr s=
k\int_{K_R}\dr\mu(t,[[M]],v)\int_0^1
\L_c\circ\psi_{-t+s}(0,M,v)\dr s=$$
$$k\int_0^1\dr s\int_{K_R}\L_c\circ\psi_{-t+s}(0,M,v)
\dr\mu(t,[[M]],v)=
k\int_{K_R}[
\L_c(t,M,v)+\a(c)
]\dr\mu(t,[[M]],v)=$$
$$k[I_c(\mu)+\a(c)]   .   \eqno (3.6)$$
The first equality above comes from the fact that $\tilde\mu$ is invariant by $\Psi$;  the second one comes from the definition of $\s^{M,v}$; the third one comes from the fact that $\tilde\mu$ is invariant by the time-one map. The fourth equality comes from the definition of $\tilde\mu$ as a push-forward, 
the fifth one is Fubini, the sixth one is the invariance of $\mu$ and the last one is the definition of $I(\mu)$.